# imuQP: An Inverse-Matrix-Updates-Based Fast QP Solver Suitable for Real-Time MPC

Victor Trường Thịnh Lâm and Mircea Lazar

*Abstract*—This paper presents a new fast active-set quadratic programming (QP) solver based on inverse matrix updates, which is suitable for real-time model predictive control (MPC). This QP solver, called *imuQP (inverse matrix update QP)*, is based on Karush-Kuhn-Tucker (KKT) conditions and is inspired by Hildreth's QP solver. An extensive convergence and optimality analysis of imuQP, including infeasibility detection, is presented. The memory footprint and computational complexity of imuQP are analyzed and compared with qpOASES, a well-known active-set QP solver in literature. Speed and accuracy of imuQP are compared with state-of-the-art active-set QP solvers by simulating a chain of spring-connected masses in MATLAB. For illustration, MPC with integral action is used, to remove offset when tracking a constant reference, with a small sampling period of $T_s = 4$ ms. Simulation results show that imuQP is suitable for fast systems — with small $T_s$ — and is competitive with state-of-the-art active-set QP solvers.

## I. Introduction

MODEL predictive control (MPC) is an advanced control technique that is very attractive for high-tech applications. However, its implementation is hampered by the need of expensive computational hardware. This is because high computational speed is necessary for solving a quadratic program (QP) online within each sampling period $T_s$ — which is usually small for fast systems — as well as high memory capacity. Recently, there were many efforts to develop MPC design and implementation methods including QP solvers for MPC that can run on inexpensive microcontrollers, which are very attractive to use on a large scale due to their low costs, but have limited memory capacity and computational speed.

Many fast QP solvers have been developed based on different methods. Commercial active-set method (ASM) QP solvers include ODYS [1], and commercial interior-point method (IPM) QP solvers include FORCES Pro, Gurobi, Mosek, CVX and CPLEX. Academic QP solvers include: OSQP [2] which is based on the alternating direction method of multipliers (ADMM), qpDUNES [3]–[5] which is based on Newton's method and combines both IPM and ASM, where affine inequality constraints are not yet supported but only box constraints [6], ASM QP solvers qpOASES [7]–[10], Hildreth's QP solver (which we will refer to as *Hildreth*) [11, ch. 2.4.4-2.4.6, p. 63-69], [12], [13], extended Hildreth [14] and IPM

V.T.T. Lâm (Electromechanics and Power Electronics group) and M. Lazar (Control Systems group) are both with the Department of Electrical Engineering, Eindhoven University of Technology, Eindhoven, the Netherlands, e-mail: `v.t.t.lam@tue.nl` and `m.lazar@tue.nl`.

QP solvers include CLP and OOQP. A heavy-ball projected primal-dual method (HBPPDM) [15] was proposed for solving sparse QPs with only box constraints on the decision variables, containing future states and control inputs. Other commercial QP solvers are provided by MATLAB, such as quadprog (IPM or ASM), mpcActiveSetSolver which is based on the active-set QP solver QPKWIK from [16], and mpcInteriorPointSolver.

In [17], an Arduino Uno has been used for input-output communication to the plant, while online MPC problems formulated as QPs (MPC-QPs) are solved on a PC with an i7 2.40-GHz Intel Core processor using MATLAB. Both MPC-QP solving and input-output communication were performed on a PC using MATLAB and SIMULINK without additional microcontrollers in [18]. Besides fast QP solvers, effort has been put into low-complexity MPC implementations on inexpensive microcontrollers. A fast gradient method for the special case of box inequality constraints, which are easier to solve than affine inequality constraints, was developed and utilized in MPC on a LEGO NXT Brick [19]. Offline MPC such as unconstrained MPC on an STM32 microcontroller [20] and Teensy 4.0 platform [21] and explicit MPC based on convex liftings on an STM32 microcontroller [22], which is a viable alternative to online MPC, were implemented. In [23], MPC was run on an Arduino for slow applications such as process control, with large sampling periods $T_s$ of 1 s and 3 s.

In Table I, the main results of real-time MPC implementations employing QP solvers are shown. In [24], qpOASES with warm start has been implemented on a PC, with code size 835 kB and average and maximum CPU time of 2.5 ms and 10 ms, respectively, for prediction horizon $N_p = 12$ and control horizon $N_c = 5$. Programmable logic controllers (PLCs) are used to run qpOASES in [25], [26] and both qpOASES and Hildreth in [27]–[29]. In [30], [31], qpOASES has been executed on an NI myRIO-1900. Hardware-in-the-loop (HIL) simulations were performed on an Atmel ARM Cortex-M3 using qpOASES and qpDUNES in [32], where sometimes the memory footprint was too big and overflowed. However, these implementations considered slow applications with large $T_s$, QP solvers with big memory footprints, expensive platforms with high computational speed and memory capacity, or low-complexity QPs with small prediction and control horizons.

This motivates the design of a new fast QP solver with low computational complexity and memory footprint, which is inspired by Hildreth due to its attractive features: *(i)* it permits simple implementation in C/C++ code and is not memory demanding, *(ii)* it can be executed fast and efficiently on inexpensive microcontrollers due to its simple operations, *(iii)* its open-source MATLAB code is publicly accessible in [11, ch. 2.4.4-2.4.6, p. 63-69]. In [33], a faster version of Hildreth called *Hildreth+'* has been designed and simulated in MATLAB on a PC. However, running Hildreth+' on an Arduino Due [34] with a small $T_s$ was still computationally challenging.

Motivated by the above, we present a new fast active-set QP solver called *imuQP (inverse matrix update QP)*. The contributions of this paper consist of *(i)* a new QP solver based on Karush-Kuhn-Tucker (KKT) conditions, with an extensive convergence and optimality analysis including infeasibility detection, *(ii)* computational complexity and memory footprint analysis of imuQP and detailed comparison with qpOASES — one of the fastest and most accurate competitive academic active-set QP solvers freely available [7], [8] — and *(iii)* comparing the CPU time and accuracy of imuQP with state-of-the-art active-set QP solvers, where MATLAB's quadprog's solution is regarded as optimal. A chain of six spring-connected masses, which is a fast system requiring a small sampling period of $T_s = 4$ ms, is simulated in MATLAB on a PC.

The remainder of this paper is structured as follows. Brief preliminaries about linear MPC and its QP formulation as well as a brief introduction of the two benchmark active-set QP solvers Hildreth and qpOASES are given in Section II. In Section III, we present the new QP solver imuQP, followed by an extensive analysis and comparison with qpOASES in Section IV. In Section V, we show simulation results and compare speed and accuracy of imuQP with state-of-the-art QP solvers. Conclusions are drawn in Section VI.

## II. PRELIMINARIES

FOR robustness, we employ integral-action linear MPC to remove offsets when tracking a constant reference. However, the developed QP solver imuQP can be used to solve any linear MPC problem or as a QP solver in sequential quadratic programming (SQP) solvers for non-linear MPC.

### A. Integral-action linear MPC with affine term

The continuous-time state-space model $\Sigma_c$ is given as

$$\Sigma_c: \begin{cases} \dot{x}_p(t) = A_c x_p(t) + \begin{bmatrix} B_c & w_c \end{bmatrix} \underbrace{\begin{bmatrix} u(t) \\ 1 \end{bmatrix}}_{\bar{u}(t)} \\ y(t) = C_c x_p(t) + \begin{bmatrix} D_c & 0 \end{bmatrix} \bar{u}(t), \end{cases} \quad (1)$$

with the plant state $x_p(t)$, the control input $u(t)$, the measured output $y(t)$ and the zero matrix 0 is of appropriate dimensions. The affine term $w_c$ usually occurs when linearizing a non-linear system around a non-zero operating point. Therefore, we need to augment $u(t)$ with an additional 1, so we can convert $\Sigma_c$ to the discrete-time state-space model $\Sigma_d$ with time instant $k$, using the MATLAB function c2d() and sampling period $T_s$

$$\Sigma_d: \begin{cases} x_p(k+1) = A_d x_p(k) + \begin{bmatrix} B_d & w_d \end{bmatrix} \underbrace{\begin{bmatrix} u(k) \\ 1 \end{bmatrix}}_{\bar{u}(k)} \\ y(k) = C_d x_p(k) + \begin{bmatrix} D_d & 0 \end{bmatrix} \bar{u}(k), \end{cases} \quad (2)$$

where we assume $u(k)$ cannot influence $y(k)$ directly, i.e. $D_c = D_d = 0$. For integral-action MPC, we define

$$\Delta x_p(k) = x_p(k) - x_p(k-1), \; \Delta u(k) = u(k) - u(k-1),$$
$$\Delta y(k+1) = y(k+1) - y(k).$$

The augmented discrete-time state-space model is given as

$$x(k+1) = \underbrace{\begin{bmatrix} A_d & 0 \\ C_d A_d & I \end{bmatrix}}_{A} x(k) + \underbrace{\begin{bmatrix} B_d \\ C_d B_d \end{bmatrix}}_{B} \Delta u(k)$$
$$y(k) = \underbrace{\begin{bmatrix} 0 & I \end{bmatrix}}_{C} x(k), \text{ where } x(k) = \begin{bmatrix} \Delta x_p(k) \\ y(k) \end{bmatrix}, \quad (3)$$

TABLE I
MEMORY FOOTPRINT AND CPU TIME OF STATE-OF-THE-ART QP SOLVERS ON VARIOUS HARDWARE. THE ABBREVIATIONS AND SYMBOLS ARE: REF.: REFERENCE; $N_p$ AND $N_c$: PREDICTION AND CONTROL HORIZON; $n_{x_p}$: STATE DIMENSION OF PLANT MODEL; PS: PLC PROGRAM/CODE SIZE; DS: PLC DATA SIZE; MAX.: MAXIMUM; AVG.: AVERAGE; DPFP: DOUBLE PRECISION FLOATING POINT; WS: WARM START; NV: NON-VOLATILE.

| Ref. | Hardware | CPU speed | Memory capacity | QP solver | $T_s$ | Memory footprint | CPU time | $N_p, N_c$ | $n_{x_p}$ |
|---|---|---|---|---|---|---|---|---|---|
| [25] | ABB AC500 PM592-ETH (PLC) | 400 MHz | User program: 4 MB RAM, User data: 4 MB integrated | qpOASES | – | PS: 372 kB, DS: 119 kB | Max. 12 ms, avg. 3 ms | 17 variables, 83 inequalities | – |
| [26] | | | | | 1 s | DPFP: PS: 1.59 MB, DS: 0.29 MB | WS + DPFP: max. 13.1 ms, avg. 4.9 ms | 24 variables, 96 inequalities, 24 bounds | – |
| [26] | | | | | | DPFP: PS: 0.37 MB, DS: 0.11 MB | WS + DPFP: max. 7.3 ms, avg. 2.0 ms | 17 variables, 78 inequalities, 17 bounds | – |
| [27] | Siemens Simatic S7-300 CPU 319-3 PN/DP (PLC) | 250 MHz | Work memory: 2048 kB, Load memory card (plug-in): 8 MB | qpOASES | 1 s | – | Max. 125 ms | $N_p = 22$, $N_c = 7$ | 4 |
| | | | | Hildreth | | | Max. 23 ms | | |
| [28] | | | | qpOASES | | 112 kB | Max. 206 ms | | |
| | | | | Hildreth | | 7 kB | Max. 51 ms | | |
| [29] | | | | qpOASES | 60 s | 51420 B | Max. 2026 ms, avg. 1308.8 ms | $N_p = 50$, $N_c = 10$ | 13 |
| | | | | Hildreth | | 53124 B | Max. 556 ms, avg. 212.4 ms | | |
| [30] | NI myRIO-1900 | 667 MHz | NV: 512 MB, DDR3: 256 MB | qpOASES | 1.430 s | – | 1430 ms | $N_p = N_c = 10$ | 4 |
| [31] | | | | | | 756 KB | | | |
| | | | | | | 707 KB | 400 ms | $N_p = N_c = 5$ | |
| [32] | Atmel ARM Cortex-M3 | 84 MHz | Flash memory: 512 kB, SRAM: 96 kB | qpOASES | 500 ms | 352 kB | Avg. 242 $\mu$s | $N_p = N_c = 5$ | 2 |
| | | | | qpDUNES | | 298 kB | Avg. 305 $\mu$s | | 6 |
| | | | | qpOASES, qpDUNES | | Memory overflow | Not available | | 9 |

where the zero matrix $0$ and identity matrix $I$ are of appropriate dimensions. Furthermore, $x(k) \in \mathbb{R}^{n_x}$, $u(k) \in \mathbb{R}^{n_u}$, $\Delta u(k) \in \mathbb{R}^{n_u}$ and $y(k) \in \mathbb{R}^{n_y}$. The matrix pair $(A, B)$ should be controllable and the matrix pair $(A, C)$ should be observable. Note that the affine term $w_d$ cancels out in (3).

Next, we define predicted inputs $\Delta U_k \in \mathbb{R}^{Nn_u}$, predicted outputs $Y_k \in \mathbb{R}^{Nn_y}$ and future reference $\mathcal{R}_k \in \mathbb{R}^{Nn_y}$ as

$$\Delta U_k = \begin{bmatrix} \Delta u_{0|k} \\ \vdots \\ \Delta u_{N-1|k} \end{bmatrix}, Y_k = \begin{bmatrix} y_{1|k} \\ \vdots \\ y_{N|k} \end{bmatrix}, \mathcal{R}_k = \begin{bmatrix} r_{1|k} \\ \vdots \\ r_{N|k} \end{bmatrix},$$

with the following relation

$$Y_k = \underbrace{\begin{bmatrix} CA \\ \vdots \\ CA^N \end{bmatrix}}_{\Phi} x(k) + \underbrace{\begin{bmatrix} CB & \cdots & 0 \\ \vdots & \ddots & \vdots \\ CA^{N-1}B & \cdots & CB \end{bmatrix}}_{\Gamma} \Delta U_k, \quad (4)$$

where $x(k) = x_{0|k}$ is the measured augmented state and $N$ is the prediction horizon, i.e. the number of samples ahead we look in the future. In this paper, the control horizon is equal to the prediction horizon, i.e. $N_c = N$. The weighting matrices $Q, R \succ 0$ penalize the stage tracking error and stage control effort, respectively, $P \succ Q$ penalizes the terminal tracking error and the future tracking error is defined as $e_{i|k} = y_{i|k} - r_{i|k}$ for $i = 0, \ldots, N$. We define the MPC cost function as

$$J(x(k), \Delta U_k) = e_{N|k}^T P e_{N|k} + \sum_{i=0}^{N-1} \left[ e_{i|k}^T Q e_{i|k} + \Delta u_{i|k}^T R \Delta u_{i|k} \right]$$
$$= e_{0|k}^T Q e_{0|k} + \mathcal{E}_k^T \Omega \mathcal{E}_k + \Delta U_k^T \Psi \Delta U_k, \quad (5)$$

where

$$\mathcal{E}_k = \begin{bmatrix} e_{1|k}^T & \cdots & e_{N|k}^T \end{bmatrix}^T = Y_k - \mathcal{R}_k,$$
$$\Omega = \mathrm{diag}(Q, \ldots, Q, P), \quad \Psi = \mathrm{diag}(R, R, \ldots, R),$$

which drives both $e_{i|k}$ and $\Delta u_{i|k}$ to $0$ when the steady-state has been reached. Substituting $Y_k$ in (4) into (5) and removing the constant terms, which does not change the minimizer, gives

$$\bar{J}(x(k), \Delta U_k) = \frac{1}{2} \Delta U_k^T E \Delta U_k + \Delta U_k^T F, \quad (6)$$

where, since $\Omega \succ 0$ and $\Psi \succ 0$, we have that

$$E = 2(\Psi + \Gamma^T \Omega \Gamma) \succ 0, \quad F = 2\Gamma^T \Omega (\Phi x(k) - \mathcal{R}_k). \quad (7)$$

Both $\bar{J}$ and $J$ have the same unique global minimizer (both are strictly convex since $E \succ 0$), so it suffices to minimize $\bar{J}$. In practice, there are always actuator constraints on the control input or safety constraints on the output. In general, constraints could also be imposed on the state $x(k)$, but we only consider

$$\Delta u_{\min} \leq \Delta u_{i|k} \leq \Delta u_{\max}, \quad u_{\min} \leq u_{i|k} \leq u_{\max},$$
$$y_{\min} \leq y_{i|k} \leq y_{\max}, \quad (8)$$

which is a common constraint formulation in MPC applications. Note that constraints in (8) translate into affine inequality constraints in the decision variables $\Delta u_{i|k}$, i.e. polytopic constraint sets, which imuQP can easily deal with.

Next, by using standard techniques (see Appendix A for the complete details), we can formulate the constrained MPC problem as a convex QP problem with the convex cost function $\bar{J}$ in (6) and affine inequality constraints in (A.5), i.e.

$$\min_{\Delta U_k} \bar{J}(x(k), \Delta U_k) = \frac{1}{2} \Delta U_k^T E \Delta U_k + \Delta U_k^T F$$
$$\text{s.t.} \quad \mathcal{L} \Delta U_k \leq d + W x(k) + V u(k-1). \quad (9)$$

*Remark 1.* If $x(k)$ is not measurable, an observer is typically used for estimating it from $y(k)$, with its dynamics

$$\hat{x}(k+1) = (A - LC)\hat{x}(k) + B\Delta u(k) + Ly(k), \quad (10)$$

where $\hat{x}(k)$ is the estimated state. Defining the estimation error $\varepsilon(k) = x(k) - \hat{x}(k)$ and its dynamics

$$\varepsilon(k+1) = (A - LC)\varepsilon(k),$$

the observer gain matrix $L$ is chosen such that the eigenvalues $|\lambda(A - LC)| < 1$ and $\varepsilon(k)$ converges to $0$ as fast as possible.

Finally, the receding horizon principle is applied as follows. First, $x(k)$ is measured (or estimated from $y(k)$ using the observer model (10)). Second, the MPC-QP (9) is solved over the prediction horizon $N$ to obtain the minimizer $\Delta U_k^*$. Finally, from $\Delta U_k^*$ the first element $\Delta u(k) = \Delta u_{0|k}^*$ is used for computing the next control input as $u(k) = \Delta u(k) + u(k-1)$. At the next time instant $k+1$, the prediction horizon $N$ is shifted forward accordingly and this process is repeated.

The MPC-QP (9) is written in the general QP form

$$\min_{\theta} \bar{J}(\theta) = \frac{1}{2} \theta^T E \theta + \theta^T F, \quad \text{s.t.} \quad M\theta \leq \gamma, \quad (11)$$

where $\theta \in \mathbb{R}^{Nn_u}$ is the decision variable and we used

$$\theta = \Delta U_k, \quad M = \mathcal{L}, \quad \gamma = d + Wx(k) + Vu(k-1). \quad (12)$$

The cost function $\bar{J}(\theta) \in \mathbb{R}$ is quadratic in $\theta$ and the inequality constraints are affine in $\theta$, which include box constraints. For consistency, we use the same letters for the matrices $E, M$ and vectors $F, \gamma$ as in [11, eq. (2.27)-(2.28), p. 53]. The Hessian matrix $E \in \mathbb{R}^{Nn_u \times Nn_u}$ and constraint matrix $M \in \mathbb{R}^{p \times Nn_u}$ are constant and are computed offline, while gradient vector $F \in \mathbb{R}^{Nn_u}$ and constraint vector $\gamma \in \mathbb{R}^p$ depend on the state $x(k)$, future reference $\mathcal{R}_k$ and previous control input $u(k-1)$.

In what follows, we denote the feasible set as $\mathcal{X} = \{\theta \mid M\theta \leq \gamma\}$ and the $n$-th constraint as $g_n$: $m_n\theta \leq \gamma_n$, where $m_n$ is the $n$-th row in $M$ and $\gamma_n$ is the $n$-th element in $\gamma$. The total number of constraints is denoted as $p$, the active set is denoted by $\mathbb{A}$ which consists of indices of active constraints and the number of active constraints is denoted as $c$. We impose upper and lower bounds on the stage variables $y_{i|k}, u_{i|k}$ and $\Delta u_{i|k}$ for $i = 0, \ldots, N-1$ (giving $2Nn_y + 4Nn_u$ constraints) and upper and lower bounds on the terminal variable $y_{N|k}$ (giving $2n_y$ constraints), which gives

$$p = 2(N+1)n_y + 4Nn_u, \quad w = \frac{p}{2}, \quad c = |\mathbb{A}| \leq w, \quad (13)$$

where $|\cdot|$ denotes the cardinality of a set, i.e. the number of elements in a set. Note that $c$ is upper bounded by $w$ since at most half of the total number of constraints can be active, i.e. are satisfied with equality. This is because the constraints consist of an equal number of lower and upper bounds which cannot be active simultaneously.

## B. Benchmark QP solvers: Hildreth and qpOASES

Two benchmark active-set QP solvers Hildreth and qpOASES, which are based on KKT conditions, will be briefly introduced, as imuQP was inspired by Hildreth because of its simplicity. The computational complexity and memory footprint of imuQP will be compared with qpOASES, one of the fastest and most accurate QP solvers freely available.

*1) Hildreth:* Mainly simple mathematical operations, divisions by a scalar and only little matrix-vector operations are used. Inputs to the solver are the matrices $E, M$ and vectors $F, \gamma$ from the QP (11) and parametrization is not required. Elements of the Lagrange multiplier $\lambda$ are updated scalar-wise until they converge numerically within some tolerance and the optimal active set $\mathbb{A}$ is then found as the indices of the non-zero components in the found $\lambda^*$. The constrained solution $\theta^*$ is obtained by using $\lambda^*$ to correct the unconstrained solution $\theta_0$. For its algorithm and implementation details, see [11, ch. 2.4.4-2.4.6, p. 63-69]. Here, the tolerance is set to $\delta = 1 \cdot 10^{-7}$, the maximum number of iterations is set to $\bar{m} = 38$.

*2) qpOASES:* Mainly matrix-vector operations and matrix updates with Cholesky and TQ factorizations and Givens plane rotations are used. Inputs are the matrices $E, M$ and vectors $F, \gamma$ from the QP (11) and $F$ and $\gamma$ are parametrized as affine functions of parameter $w$ (i.e. it consists of the state $x(k)$, future reference $\mathcal{R}_k$ and previous control input $u(k-1)$). The Lagrange multiplier $\lambda$ and constrained solution $\theta$ are computed exactly based on the current active set $\mathbb{A}$. By adding and removing constraints, the active set $\mathbb{A}$ converges to an optimal one. For its algorithm and details, see [7]–[10].

## III. IMUQP ACTIVE-SET QP SOLVER

IN this section, we present the designed QP solver called imuQP. We start with the KKT conditions, which are the key equations in imuQP. Then, we will present the algorithm of imuQP and discuss *additive inverse matrix updates (AIMU)* and *subtractive inverse matrix updates (SIMU)*, which are some of the key features that make imuQP very fast.

For the QP in (11), let us assume that the optimal constrained solution (primal solution) $\theta^*$ and optimal Lagrange multiplier (dual solution) $\lambda^* \in \mathbb{R}_{\geq 0}^p$ were found. The optimal active set is also assumed to be known, i.e. these constraints are satisfied with equality. Moreover, let the indices of active and inactive constraints be denoted by $\mathbb{A}$ and $\mathbb{I} = \{1, \ldots, p\} \setminus \mathbb{A}$, respectively. Then from [11, eq. (2.40)-(2.43), p. 58]

$$M_{\mathbb{A}} \theta^* - \gamma_{\mathbb{A}} = 0, \quad \lambda_{\mathbb{A}}^* > 0, \qquad (14a)$$
$$M_{\mathbb{I}} \theta^* - \gamma_{\mathbb{I}} < 0, \quad \lambda_{\mathbb{I}}^* = 0, \qquad (14b)$$

where the symbols $\mathbb{A}$ and $\mathbb{I}$ are adopted from [7], [8]. Then, (14) is equivalent to

$$M\theta^* - \gamma \leq 0, \quad (\lambda^*)^T (M\theta^* - \gamma) = 0, \quad \lambda^* \geq 0, \quad (15)$$

where $M_{\mathbb{A}} \in \mathbb{R}^{c \times Nn_u}, \gamma_{\mathbb{A}} \in \mathbb{R}^c$ and $\lambda_{\mathbb{A}}^* \in \mathbb{R}_{\geq 0}^c$ (or $M_{\mathbb{I}}, \gamma_{\mathbb{I}}$ and $\lambda_{\mathbb{I}}^*$) contain the $n$-th rows in $M$ and the $n$-th elements in $\gamma$ and $\lambda^*$, respectively, with $n \in \mathbb{A}$ (or $n \in \mathbb{I}$).

By leaving out the inactive constraints, this boils down to solving the equality-constrained quadratic program (EQP)

$$\min_{\theta} \bar{J}(\theta) = \frac{1}{2} \theta^T E \theta + \theta^T F, \quad \text{s.t. } M_{\mathbb{A}} \theta = \gamma_{\mathbb{A}}.$$

From [11, eq. (2.30), p. 54], the Lagrangian function for this EQP is constructed as

$$\mathcal{L}(\theta, \lambda_{\mathbb{A}}) = \frac{1}{2} \theta^T E \theta + \theta^T F + \lambda_{\mathbb{A}}^T (M_{\mathbb{A}} \theta - \gamma_{\mathbb{A}}),$$

and the following equations are satisfied [11, eq. (2.31)-(2.32), p. 55]

$$\frac{\partial \mathcal{L}(\theta^*, \lambda_{\mathbb{A}}^*)}{\partial \theta} = E\theta^* + F + M_{\mathbb{A}}^T \lambda_{\mathbb{A}}^* = 0 \qquad (16a)$$
$$\frac{\partial \mathcal{L}(\theta^*, \lambda_{\mathbb{A}}^*)}{\partial \lambda} = M_{\mathbb{A}} \theta^* - \gamma_{\mathbb{A}} = 0, \qquad (16b)$$

and therefore can be expressed as [11, eq. (2.33)-(2.34), p. 55]

$$\lambda_{\mathbb{A}}^* = -H_{\mathbb{A}}^{-1} K_{\mathbb{A}}, \quad \theta^* = \theta_0 - E^{-1} M_{\mathbb{A}}^T \lambda_{\mathbb{A}}^*. \qquad (17)$$

Here,

$$H_{\mathbb{A}} = M_{\mathbb{A}} E^{-1} M_{\mathbb{A}}^T, \quad K_{\mathbb{A}} = \gamma_{\mathbb{A}} - M_{\mathbb{A}} \theta_0, \qquad (18)$$

and the unconstrained solution $\theta_0$ is obtained from

$$\nabla \bar{J}(\theta_0) = E\theta_0 + F = 0 \quad \rightarrow \quad \theta_0 = -E^{-1} F. \qquad (19)$$

The necessary KKT conditions are formed by (15) and

$$E\theta^* + F + M^T \lambda^* = 0, \qquad (20)$$

which is a generalized form of (16a) since $\lambda_{\mathbb{I}}^* = 0$, see (14b).

The imuQP algorithm for solving the general QP in (11) is presented in Algorithm 1. Here, the optimal active set $\mathbb{A}^*$ and its size $c^*$ correspond to the optimal $\theta^*$ and $\lambda_{\mathbb{A}}^*$.

---

**Algorithm 1** imuQP (inverse matrix update QP).

---

**Inputs:**
| | |
|---|---|
| $E \in \mathbb{R}^{Nn_u \times Nn_u}$ | Hessian matrix |
| $F \in \mathbb{R}^{Nn_u}$ | gradient vector |
| $M \in \mathbb{R}^{p \times Nn_u}$ | constraint matrix |
| $\gamma \in \mathbb{R}^p$ | constraint vector |
| $H \in \mathbb{R}^{p \times p}$ | dual Hessian matrix |

**Outputs:**
| | |
|---|---|
| $\theta^* \in \mathbb{R}^{Nn_u}$ | optimal constrained solution |
| $\lambda_{\mathbb{A}}^* \in \mathbb{R}_{\geq 0}^{c^*}$ | optimal Lagrange multiplier |
| $\mathbb{A}^* \in \mathbb{Z}_{\geq 0}^{c^*}$ | constraint indices of optimal active set |
| $c^* \in \mathbb{Z}_{\geq 0}$ | $|\mathbb{A}^*|$, i.e. size of optimal active set |
| $m^* \in \mathbb{Z}_{> 0}$ | number of iterations performed |

**Steps:**
*I. Initialization phase*

1. Compute unconstrained solution vector $\theta_0 \in \mathbb{R}^{Nn_u}$, i.e.

$$\theta_0 = -E^{-1} F. \qquad (21)$$

2. Compute constraint-violation vector $K_0 \in \mathbb{R}^p$, i.e.

$$K_0 = \gamma - M\theta_0. \qquad (22)$$

3. a. Set iteration counter to $m = 1$.
   b. From $K_0$, extract $K_{\min}$ with constraint index $j$, i.e.

$$K_{\min} = \min(K_0) \quad \text{and} \quad j = \arg\min(K_0). \qquad (23)$$

   c. Initialize: $\lambda_{\mathbb{A}} \in \mathbb{R}^w$ and $\mathbb{A} \in \mathbb{Z}_{\geq 0}^w$ filled with zeros, $c = 0$.
   d. **if** $K_{\min} \geq 0$ **then** set $\theta^* = \theta_0$ and <u>return</u> (no violation of constraint) **else** continue (violation of constraint) **end**
   e. Initialize: $H_{\mathbb{A}}^{-1} \in \mathbb{R}^{w \times w}$ filled with zeros.

*II. Update phase*

4. **while** $K_{\min} < 0$ (constraint violated) **then**
   a. **if** $j \in \mathbb{A}$, i.e. $\mathbb{A}$ contains $j$ **then** <u>break</u> **else**
      i. Increase iteration counter $m = m + 1$.
      ii. Set $c = c + 1$. Add constraint index $j$ to $\mathbb{A}$, i.e. $\mathbb{A} = \mathbb{A} \cup \{j\}$.
   **end**
   b. **if** $c = 1$, i.e. $|\mathbb{A}| = 1$ ($\mathbb{A}$ contains 1 index) **then**
   $$H_{\mathbb{A}}^{-1} = \frac{1}{h_{jj}} \quad \text{and} \quad \lambda_{\mathbb{A}} = -H_{\mathbb{A}}^{-1} \cdot K_{\min}, \quad (24)$$
   where $h_{jj}$ is the $(j,j)$-th element in $H$ **else**
      i. Check if constraint index $j$ depends linearly on the constraints that are present in $\mathbb{A}$
      $$y = H_{\text{old}}^{-1}\mathbf{h} \quad \text{and} \quad q = h - \mathbf{h}^T y. \quad (25)$$
      ii. **if** $q = 0$ **then** linear dependence:
         A. Extract $y_{\max}$, which is the $f$-th entry in $y$, i.e.
         $$y_{\max} = \max(y) \quad \text{and} \quad f = \arg\max(y). \quad (26)$$
         B. **if** $y_{\max} \leq 0$ **then** <u>break</u> (QP infeasible) **else** continue
         **end**
         C. Set $c = c - 1$. Remove the $f$-th entry in $\mathbb{A}$, i.e. $\mathbb{A} = \mathbb{A}\backslash\{\mathbb{A}_f\}$.
         D. Update $H_{\mathbb{A}}^{-1}, \lambda_{\mathbb{A}}$ using SIMU with $f$-th entry removed.
         E. Update $H_{\mathbb{A}}^{-1}, \lambda_{\mathbb{A}}$ using AIMU with constraint index $j$ appended.
      **else** no linear dependence: update $H_{\mathbb{A}}^{-1}, \lambda_{\mathbb{A}}$ using AIMU with constraint index $j$ appended **end**
   **end**
   c. Extract $\lambda_{\min}$, which is the $i$-th entry in $\lambda_{\mathbb{A}}$, i.e.
   $$\lambda_{\min} = \min(\lambda_{\mathbb{A}}) \quad \text{and} \quad i = \arg\min(\lambda_{\mathbb{A}}). \quad (27)$$
   d. **while** $\lambda_{\min} < 0$ (constraint becomes inactive) **then**
      i. Increase iteration counter $m = m + 1$.
      ii. Set $c = c - 1$. Remove the $i$-th entry in $\mathbb{A}$, i.e. $\mathbb{A} = \mathbb{A}\backslash\{\mathbb{A}_i\}$.
      iii. Update $H_{\mathbb{A}}^{-1}, \lambda_{\mathbb{A}}$ using SIMU with $i$-th entry removed.
      iv. Extract $\lambda_{\min}$, which is the $i$-th entry in $\lambda_{\mathbb{A}}$, i.e.
      $$\lambda_{\min} = \min(\lambda_{\mathbb{A}}) \quad \text{and} \quad i = \arg\min(\lambda_{\mathbb{A}}). \quad (28)$$
   **end**
   e. Compute new constraint-violation vector $K \in \mathbb{R}^p$, i.e.
   $$K = K_0 + H_{\mathbb{A},\text{col}}\lambda_{\mathbb{A}}, \quad (29)$$
   where $H_{\mathbb{A},\text{col}}$ consists of $n$-th columns in $H$, with $n \in \mathbb{A}$.
   f. From $K$, extract $K_{\min}$ with constraint index $j$, i.e.
   $$K_{\min} = \min(K) \quad \text{and} \quad j = \arg\min(K). \quad (30)$$
**end**

*III. Final phase*

5. Compute the optimal constrained solution $\theta^* \in \mathbb{R}^{Nn_u}$, as a correction of $\theta_0$ using the found $\mathbb{A}^*$ and $\lambda_{\mathbb{A}}^*$, i.e.
$$\theta^* = \theta_0 - E^{-1}M_{\mathbb{A}}^T\lambda_{\mathbb{A}}^*. \quad (31)$$

## A. Implementation details

Both imuQP and Hildreth are active-set QP solvers, based on KKT conditions where the optimal $\theta^*$ is computed from the optimal $\lambda^*$ and the unconstrained solution $\theta_0$. However, there are three main differences between the two QP solvers. In what follows, $H_{\mathbb{A}} \in \mathbb{R}^{c \times c}$ contains the $n$-th rows and $n$-th columns of the dual Hessian matrix $H \in \mathbb{R}^{p \times p}$, which is
$$H = ME^{-1}M^T, \quad (32)$$
and $K_{\mathbb{A}}$ and $\lambda_{\mathbb{A}} \in \mathbb{R}^c$ contain the $n$-th elements in $K_0$ in (22) and $\lambda \in \mathbb{R}^p$, respectively, with $n \in \mathbb{A}$.

First, for compact notation we write $E^{-1}$ in the expressions for $H$ in (32), $\theta_0$ in (21) and $\theta^*$ in (31), but is implemented differently in both QP solvers. The Hildreth implementation in [11, ch. 2.4.4-2.4.6, p. 63-69] uses the backslash operator (\) to solve a system of linear equations such as $Ax = b$, i.e. $x = A\backslash b$ in MATLAB and computes online: *(i)* $H = M(E\backslash M^T)$, *(ii)* $\theta_0 = -E\backslash F$ and *(iii)* $\theta^* = -E\backslash F - E\backslash(M^T\lambda^*)$. Instead, imuQP factorizes $E$ using Gaussian elimination (without pivoting) and uses forward and backward substitution − which are fast and accurate − to compute: *(i)* offline the constant matrix $H$ as $H = M\eta$, where $\eta = \begin{bmatrix} \eta_1 & \cdots & \eta_p \end{bmatrix} \in \mathbb{R}^{Nn_u \times p}$. Each $\eta_q \in \mathbb{R}^{Nn_u}$ (the $q$-th column in $\eta$) for $q = 1,\ldots,p$ is the solution to $E\eta_q = m_q^T$, where $m_q$ is the $q$-th row of $M$, *(ii)* online the solution to $E\theta_0 = -F$ and *(iii)* online $\theta^* = \theta_0 - \kappa$, where $\kappa$ is the solution to $E\kappa = M_{\mathbb{A}}^T\lambda_{\mathbb{A}}^*$.

Second, ideally we can find $\lambda_{\mathbb{A}}^*$ immediately from (17) if we knew the optimal active set $\mathbb{A}^*$ in advance [35, ch. 12.3, p. 361], which is difficult. Instead, Hildreth updates the elements of $\lambda$ scalar-wise until they converge numerically within some tolerance. Similarly to qpOASES, imuQP adds the *most violated constraint* to the active set $\mathbb{A}$ or removes the *most inactive constraint* from $\mathbb{A}$ until the solution is both feasible and optimal within some tolerance. In [36, alg. 16.3, p. 472] and [37] a similar active-set method is described, but requires an initial feasible solution, unlike imuQP. For the current active set $\mathbb{A}$, $H_{\mathbb{A}}^{-1}$ and $\lambda_{\mathbb{A}}$ are computed exactly with AIMU and SIMU, using very cheap operations. $H_{\mathbb{A}}^{-1}$, $\lambda_{\mathbb{A}}$ and $\mathbb{A}$ are initialized as zero matrices and zero vectors with fixed sizes $\mathbb{R}^{w \times w}$, $\mathbb{R}^w$ and $\mathbb{R}^w$, respectively, to avoid dynamic array sizes. Only the first $c$ columns and the first $c$ rows in $H_{\mathbb{A}}^{-1}$ and the first $c$ elements in $\lambda_{\mathbb{A}}$ and $\mathbb{A}$ are actually used.

Third, the Hildreth implementation in [11, ch. 2.4.4-2.4.6, p. 63-69] successively checks whether constraint $i$ is violated or not, i.e. $m_i\theta_0 > \gamma_i$ for $i = 1,\ldots,p$. The vector $K_0 = \gamma + M(E\backslash F) \in \mathbb{R}^p$ is later computed from scratch, again using the backslash operator. Instead, imuQP computes $K_0$ in (22) in a vectorized fashion and checks if $\min(K_0) < 0$. Note that $\theta_0$ was already computed in (21) and reused for computing $K_0$ in (22), $K_{\mathbb{A}}$ in (18) and $\theta^*$ in (31) and (17). Moreover, $K \in \mathbb{R}^p$ in (29) is computed by reusing the already computed $K_0$ in (22), $y$ is reused in $q$ in (25) and $q$ in (25) in turn is reused in the denominators in $H_{\text{new}}^{-1}$ in (34) and $\lambda_{\text{new}}$ in (36). Note that for better readability, we use the notation $\theta_0$ in steps 1, 2 and 5 of Algorithm 1. However, in the implementation the output variable $\theta^*$ is used and step 3.d just returns $\theta^*$.

The following key features of imuQP make it fast and accurate. First, infeasiblity of a QP can be easily detected *for*

*free* as it is part of the inverse matrix update process – without any additional algorithm. Second, Gaussian elimination with forward and backward substitution, AIMU, SIMU, infeasibility detection of a QP and other parts in imuQP are implemented as efficiently as possible, i.e. with the least number of arithmetic operations, reuse of previously computed variables, efficient **for** loops and vectorizations. Moreover, those operations involve only simple mathematical operations such as matrix-vector multiplications, divisions by a scalar, indexing vectors and matrices and finding the minimum or maximum value and its corresponding index in a vector.

### B. Inverse matrix updates

AIMU and SIMU are used to update $H_{\mathbb{A}}^{-1}$ and $\lambda_{\mathbb{A}}$ whenever a single constraint is added or removed from the active set $\mathbb{A}$. Inverse matrix updates are computationally cheaper than inverting a matrix from scratch because they make use of previously known matrix inverses and only involve matrix-vector multiplications and divisions by a scalar, which are simple mathematical operations. In what follows, let $\mathbb{A}_{\text{old}}^m$ and $\mathbb{A}_{\text{new}}^m$ denote the active set before and after a constraint has been added or removed at iteration $m$, respectively, and $H_{\text{old}}, K_{n \in \mathbb{A}_{\text{old}}}, \mathbf{k}_{\text{old}}$ and $\lambda_{\text{old}}$ (or $H_{\text{new}}, K_{n \in \mathbb{A}_{\text{new}}}, \mathbf{k}_{\text{new}}$ and $\lambda_{\text{new}}$) contain the $n$-th rows and $n$-th columns from $H$ in (32) and the $n$-th elements from $K$ in (29), $K_0$ in (22) and $\lambda$, respectively, with $n \in \mathbb{A}_{\text{old}}^m$ (or $n \in \mathbb{A}_{\text{new}}^m$). Note that $\mathbb{A}_{\text{old}}^{m+1} = \mathbb{A}_{\text{new}}^m$.

*1) AIMU:* Suppose for $\mathbb{A}_{\text{old}}^m$ we know $H_{\text{old}}$, its inverse $H_{\text{old}}^{-1}$ and $\mathbf{k}_{\text{old}}$. Now, constraint index $j$ is added to the active set, i.e. $\mathbb{A}_{\text{new}}^m = \mathbb{A}_{\text{old}}^m \cup \{j\}$ and we append $H_{\text{old}}$ with $\mathbf{h}$ and $h$ and append $\mathbf{k}_{\text{old}}$ with $k$ (highlighted in green) as follows

$$H_{\text{new}} = \begin{bmatrix} H_{\text{old}} & \mathbf{h} \\ \mathbf{h}^T & h \end{bmatrix}, \quad \mathbf{k}_{\text{new}} = \begin{bmatrix} \mathbf{k}_{\text{old}} \\ k \end{bmatrix}, \quad (33)$$

where $h$ is the $(j,j)$-th element from $H$ and $\mathbf{h}$ (or $\mathbf{h}^T$) consists of the $(n,j)$-th (or $(j,n)$-th) elements in $H$ with $n \in \mathbb{A}_{\text{old}}^m$, respectively, and $k$ is the $j$-th element from $K_0$. Then, from [38]–[40] and [41, ch. 1: remark, p. 2], we obtain

$$H_{\text{new}}^{-1} = \begin{bmatrix} H_{\text{old}}^{-1} + \frac{H_{\text{old}}^{-1} \mathbf{h} \mathbf{h}^T H_{\text{old}}^{-1}}{h - \mathbf{h}^T H_{\text{old}}^{-1} \mathbf{h}} & -\frac{H_{\text{old}}^{-1} \mathbf{h}}{h - \mathbf{h}^T H_{\text{old}}^{-1} \mathbf{h}} \\ -\frac{\mathbf{h}^T H_{\text{old}}^{-1}}{h - \mathbf{h}^T H_{\text{old}}^{-1} \mathbf{h}} & \frac{1}{h - \mathbf{h}^T H_{\text{old}}^{-1} \mathbf{h}} \end{bmatrix}, \quad (34)$$

where $H_{\text{new}}^{-1}$ is expressed in terms of $H_{\text{old}}^{-1}$, $\mathbf{h}$ and $h$ which are known. From (17), we can derive

$$\lambda_{\text{old}} = -H_{\text{old}}^{-1} \mathbf{k}_{\text{old}}, \quad \lambda_{\text{new}} = -H_{\text{new}}^{-1} \mathbf{k}_{\text{new}}, \quad (35)$$

and by substituting $H_{\text{new}}^{-1}$ in (34) and $\mathbf{k}_{\text{new}}$ in (33) into $\lambda_{\text{new}}$ in (35), we obtain

$$\lambda_{\text{new}} = \begin{bmatrix} \lambda_{\text{old}} + \frac{H_{\text{old}}^{-1} \mathbf{h}}{h - \mathbf{h}^T H_{\text{old}}^{-1} \mathbf{h}} \left(k + \mathbf{h}^T \lambda_{\text{old}}\right) \\ -\frac{1}{h - \mathbf{h}^T H_{\text{old}}^{-1} \mathbf{h}} \left(k + \mathbf{h}^T \lambda_{\text{old}}\right) \end{bmatrix}, \quad (36)$$

where $\lambda_{\text{new}}$ is expressed in terms of $\lambda_{\text{old}}, H_{\text{old}}^{-1}, \mathbf{h}, h$ and $k$ which are known.

*Remark 2.* Initially, when $\mathbb{A}_{\text{old}}^m = \varnothing$ and $\mathbb{A}_{\text{new}}^m = j$, then

$$H_{\text{new}}^{-1} = \frac{1}{h}, \quad \lambda_{\text{new}} = -H_{\text{new}}^{-1} k.$$

*2) SIMU:* Again, suppose for $\mathbb{A}_{\text{old}}^m$ we know $H_{\text{old}}$, its inverse $H_{\text{old}}^{-1}$ and $\mathbf{k}_{\text{old}}$. Now, constraint index $\mathbb{A}_i$ (the $i$-th entry in $\mathbb{A}_{\text{old}}^m$) is removed from the active set, i.e. $\mathbb{A}_{\text{new}}^m = \mathbb{A}_{\text{old}}^m \setminus \{\mathbb{A}_i\}$ and we remove the $i$-th row and $i$-th column from $H_{\text{old}}$ and the $i$-th element from $\mathbf{k}_{\text{old}}$. We partition $H_{\text{old}}^{-1}, \mathbf{k}_{\text{old}}$ and $\lambda_{\text{old}}$ and highlight the $i$-th row and $i$-th column in $H_{\text{old}}^{-1}$ and the $i$-th element in $\mathbf{k}_{\text{old}}$ and $\lambda_{\text{old}}$ in green as follows (note that the symbols $h$ and $k$ in (37) are different from those in (33))

$$H_{\text{old}}^{-1} = \begin{bmatrix} H_1 & \mathbf{h}_1 & H_2 \\ \mathbf{h}_1^T & h & \mathbf{h}_2^T \\ H_2^T & \mathbf{h}_2 & H_3 \end{bmatrix}, \quad (37a)$$

$$\mathbf{k}_{\text{old}} = \begin{bmatrix} \mathbf{k}_1 \\ k \\ \mathbf{k}_2 \end{bmatrix}, \quad \lambda_{\text{old}} = \begin{bmatrix} \boldsymbol{\lambda}_1 \\ \lambda \\ \boldsymbol{\lambda}_2 \end{bmatrix}. \quad (37b)$$

Then, from [38], [40]

$$H_{\text{new}}^{-1} = \begin{bmatrix} H_1 & H_2 \\ H_2^T & H_3 \end{bmatrix} - \frac{1}{h} \begin{bmatrix} \mathbf{h}_1 \\ \mathbf{h}_2 \end{bmatrix} \begin{bmatrix} \mathbf{h}_1^T & \mathbf{h}_2^T \end{bmatrix}, \quad (38)$$

where $H_{\text{new}}^{-1}$ is expressed in terms of $H_1, H_2, H_3, \mathbf{h}_1, \mathbf{h}_2$ and $h$ which are known. Substituting $H_{\text{old}}^{-1}, \mathbf{k}_{\text{old}}$ and $\lambda_{\text{old}}$ in (37) into $\lambda_{\text{old}}$ in (35) and using $\lambda_{\text{new}}$ in (35) with $H_{\text{new}}^{-1}$ in (38) gives

$$\lambda_{\text{new}} = \begin{bmatrix} \boldsymbol{\lambda}_1 \\ \boldsymbol{\lambda}_2 \end{bmatrix} - \frac{\lambda}{h} \begin{bmatrix} \mathbf{h}_1 \\ \mathbf{h}_2 \end{bmatrix}, \quad \text{with } \mathbf{k}_{\text{new}} = \begin{bmatrix} \mathbf{k}_1 \\ \mathbf{k}_2 \end{bmatrix}, \quad (39)$$

where $\lambda_{\text{new}}$ is expressed in terms of $\boldsymbol{\lambda}_1, \boldsymbol{\lambda}_2, \lambda, h, \mathbf{h}_1$ and $\mathbf{h}_2$ which are known. Note that $\lambda = -\mathbf{h}_1^T \mathbf{k}_1 - hk - \mathbf{h}_2^T \mathbf{k}_2$.

## IV. CONVERGENCE AND COMPLEXITY ANALYSIS

**W**E will show that imuQP is guaranteed to converge to the optimal solution in three parts. First, we assume the QP in (11) is feasible and will show that imuQP terminates within a finite number of iterations, i.e. cycling does not occur. Second, we will show that the constrained solution $\theta^* \in \mathcal{X}$ (feasibility) and that the Lagrange multiplier $\lambda_{\mathbb{A}}^* \geq 0$ (optimality) for the optimal active set $\mathbb{A}^*$. Third, we will show how to check linear dependence of constraints and infeasibility of a QP, which is part of the inverse matrix update process.

*1) Non-cycling of iterations:* Let us assume the following.

**Assumption 1.1.** *The QP in (11) is feasible, i.e. $\mathcal{X} \neq \varnothing$ and there is a solution $\theta \in \mathcal{X}$.*

**Assumption 1.2.** *When adding a new constraint $j$ to $\mathbb{A}_{\text{old}}^m$ such that $\mathbb{A}_{\text{new}}^m = \mathbb{A}_{\text{old}}^m \cup \{j\}$, step 4.b in Algorithm 1 always ensures that $M_{\mathbb{A}}$ has full row rank, i.e. it has linearly independent rows. Since $E \succ 0$, then $H_{\mathbb{A}} = M_{\mathbb{A}} E^{-1} M_{\mathbb{A}}^T \succ 0$ [42, obs. 7.1.8(b), p. 431].*

The upcoming three lemmas show that adding or removing a constraint cannot result in adding or removing the same constraint again in the next iteration.

**Lemma 1.1.** *Suppose at iteration $m$, either (i) constraint index $j$ is added to $\mathbb{A}_{\text{old}}^m$ such that $\mathbb{A}_{\text{new}}^m = \mathbb{A}_{\text{old}}^m \cup \{j\}$, or (ii) a constraint has been removed from $\mathbb{A}_{\text{old}}^m$ such that $\mathbb{A}_{\text{new}}^m = j$ contains only one constraint. Then it holds that $\lambda_c > 0$, where $\lambda_c$ is the last element in $\lambda_{\mathbb{A}} \in \mathbb{R}^c$, and at iteration $m+1$ constraint index $j$ cannot be removed from $\mathbb{A}_{\text{old}}^{m+1} = \mathbb{A}_{\text{new}}^m$.*

*Proof.* (i) Case 1: suppose at iteration $m = 1$ we have $\mathbb{A}_{\text{new}}^1 = \varnothing$ and the (unconstrained) solution $\theta_1 = \theta_0$, which violates the $j$-th constraint the most, i.e. $m_j \theta_1 > \gamma_j$. At iteration $m = 2$, we have $\mathbb{A}_{\text{new}}^2 = \mathbb{A}_{\text{old}}^2 \cup \{j\} = j$ (to be shown in Lemma 2.2.(i)) and the solution $\theta_2$. Since $|\mathbb{A}_{\text{new}}^2| = 1$, from (24) we have

$$\lambda_c = \lambda_\mathbb{A} = -\frac{K_{\min}}{h_{jj}} > 0, \qquad (40)$$

where $K_{\min} < 0$ in (23) and $h_{jj} > 0$ (because $H_\mathbb{A} \succ 0$, see Assumption 1.2).

Case 2: suppose at iteration $m \geq 2$ we have $\mathbb{A}_{\text{new}}^m$ and the solution $\theta_m$, which violates the $j$-th constraint the most, which is linearly independent on the $n$-th constraints with $n \in \mathbb{A}_{\text{new}}^m$ (from Assumption 1.2). At iteration $m+1$, we have $\mathbb{A}_{\text{new}}^{m+1} = \mathbb{A}_{\text{old}}^{m+1} \cup \{j\}$ (to be shown in Lemma 2.2.(i)) and the solution $\theta_{m+1}$. Since $|\mathbb{A}_{\text{new}}^{m+1}| > 1$, from the last component in (36)

$$\lambda_c = -\frac{1}{h - \mathbf{h}^T H_{\text{old}}^{-1} \mathbf{h}} \left( k + \mathbf{h}^T \lambda_{\text{old}} \right) > 0,$$

where $k + \mathbf{h}^T \lambda_{\text{old}} = K_{\min} < 0$ in (30) from iteration $m$ and $h - \mathbf{h}^T H_{\text{old}}^{-1} \mathbf{h} > 0$ (to be shown in Lemma 3.1.(ii)).

*(ii)* Suppose the opposite is true, i.e. $\lambda_c < 0$, which means we can relax the $j$-th constraint and move to $\bar{\theta}_m = \theta_m + \delta\theta$ such that $\bar{J}(\bar{\theta}_m) < \bar{J}(\theta_m)$ with $\bar{J}(\theta)$ in (11) while $\bar{\theta}_m \in \mathcal{X}$ (as will be shown in Lemma 2.3.(i)). But since $\mathbb{A}_{\text{new}}^{m+1} = \mathbb{A}_{\text{old}}^{m+1}\setminus\{j\} = \varnothing$, it means the unconstrained solution $\theta_0$ in (21) should also satisfy the $j$-th constraint. This is a contradiction, because imuQP would have terminated as soon as $\theta_0$ satisfies all constraints at step 3.d in Algorithm 1.

Since $\lambda_c \not< 0$, the **while** condition in step 4.d in Algorithm 1 ensures constraint index $j$ is not removed from $\mathbb{A}_{\text{old}}^{m+1} = \mathbb{A}_{\text{new}}^m$ at iteration $m+1$. ∎

**Lemma 1.2.** *Given $\mathbb{A}_{\text{old}}^m$ at iteration $m$, then (i) it holds that $m_n \theta_m = \gamma_n$ with $n \in \mathbb{A}_{\text{old}}^m$ and solution $\theta_m$, and (ii) if constraint index $j$ was already added to $\mathbb{A}_{\text{old}}^m$, then at iteration $m+1$ it cannot be added again to $\mathbb{A}_{\text{old}}^{m+1} = \mathbb{A}_{\text{new}}^m$.*

*Proof.* (i) The constraint-violation vector $K$ in (29) is derived by substituting $\theta^*$ in (31) into

$$K = \gamma - M\theta_m = \underbrace{\gamma - M\theta_0}_{K_0} + \underbrace{ME^{-1}M_\mathbb{A}^T}_{H_{\mathbb{A},\text{col}}} \lambda_\mathbb{A}, \qquad (41)$$

with $K_0$ in (22) and $\lambda_\mathbb{A}^*$ and $\theta^*$ are replaced with $\lambda_\mathbb{A}$ and $\theta_m$, respectively. With $\lambda_{\text{old}}$ in (35), we can easily see that

$$K_{n \in \mathbb{A}_{\text{old}}} = \mathbf{k}_{\text{old}} + H_{\text{old}} \lambda_{\text{old}} = \mathbf{k}_{\text{old}} - H_{\text{old}} H_{\text{old}}^{-1} \mathbf{k}_{\text{old}} = 0. \quad (42)$$

From $K$ in (41), we can conclude that

$$K_{n \in \mathbb{A}_{\text{old}}} = \gamma_n - m_n \theta_m = 0 \quad \text{or} \quad m_n \theta_m = \gamma_n.$$

*(ii)* Case 1 (from Lemma 1.1.(i)): from $K$ in (29) and $\lambda_\mathbb{A}$ in (40) we have that

$$K_{n \in \mathbb{A}_{\text{new}}} = K_j = K_{\min} + h_{jj} \lambda_\mathbb{A} = K_{\min} - h_{jj} \cdot \frac{K_{\min}}{h_{jj}} = 0,$$

where $K_j$ and $K_{\min}$ in (23) are the $j$-th element in $K$ in (29) and $K_0$ in (22), respectively.

Case 2 (from Lemma 1.1.(i)): At iteration $m$, from $K$ in (29) and using $H_{\text{new}}$ and $\mathbf{k}_{\text{new}}$ in (33) and $\lambda_{\text{new}}$ in (36)

$$K_{n \in \mathbb{A}_{\text{new}}} = \begin{bmatrix} K_{n \in \mathbb{A}_{\text{old}}} \\ K_j \end{bmatrix} = \mathbf{k}_{\text{new}} + H_{\text{new}} \lambda_{\text{new}}$$
$$= \begin{bmatrix} \mathbf{k}_{\text{old}} + H_{\text{old}} \lambda_{\text{old}} \\ k - k \end{bmatrix} = \begin{bmatrix} 0 \\ 0 \end{bmatrix}, \qquad (43)$$

where $K_j = K_{\min} < 0$, $K_{\min}$ is from (30) and we used that $\mathbf{k}_{\text{old}} + H_{\text{old}} \lambda_{\text{old}} = 0$ from (42). We obtain the same result by using $\lambda_{\text{new}}$ in (35) to see that

$$K_{n \in \mathbb{A}_{\text{new}}} = \mathbf{k}_{\text{new}} + H_{\text{new}} \lambda_{\text{new}} = \mathbf{k}_{\text{new}} - H_{\text{new}} H_{\text{new}}^{-1} \mathbf{k}_{\text{new}} = 0.$$

Since $K_j \not< 0$, the **while** condition in step 4 of Algorithm 1 ensures constraint index $j$ is not added again to $\mathbb{A}_{\text{old}}^{m+1} = \mathbb{A}_{\text{new}}^m$ at iteration $m+1$. The check in step 4.a in Algorithm 1 is required for numerical precision reasons, i.e. the right-hand side of (43) might be a very small negative number and may result in adding a constraint $j$ that is already in $\mathbb{A}_{\text{new}}^m$. ∎

**Lemma 1.3.** *Suppose at iteration $m \geq 4$, constraint index $\mathbb{A}_i$ (the $i$-th entry in $\mathbb{A}_{\text{old}}^m$) is removed from $\mathbb{A}_{\text{old}}^m$ such that $\mathbb{A}_{\text{new}}^m = \mathbb{A}_{\text{old}}^m \setminus \{\mathbb{A}_i\}$. Then it cannot be added again to $\mathbb{A}_{\text{old}}^{m+1} = \mathbb{A}_{\text{new}}^m$.*

*Proof.* From Lemma 1.1.(i), constraints cannot be removed from $\mathbb{A}_{\text{old}}^m$ at iterations $m = 1, 2, 3$. Therefore, we consider iteration $m \geq 4$ with $|\mathbb{A}_{\text{old}}^m| > 1$, where we partition $H_{\text{old}}$ and highlight its $i$-th row and $i$-th column in green as follows

$$H_{\text{old}} = \begin{bmatrix} B_1 & \mathbf{b}_1 & B_2 \\ \mathbf{b}_1^T & b & \mathbf{b}_2^T \\ B_2^T & \mathbf{b}_2 & B_3 \end{bmatrix}. \qquad (44)$$

The $\mathbb{A}_i$-th row in (29) gives

$$K_{\mathbb{A}_i} = k + \mathbf{b}_1^T \boldsymbol{\lambda}_1 + b\lambda + \mathbf{b}_2^T \boldsymbol{\lambda}_2 = 0, \qquad (45)$$

where the dot product of the $i$-th row of $H_{\text{old}}$ in (44) with $\lambda_{\text{old}}$ in (37b) is taken, $K_{\mathbb{A}_i}$ denotes the $\mathbb{A}_i$-th element of $K$ in (29) and $k$ denotes the $\mathbb{A}_i$-th element of $K_0$ in (22). We know that $H_{\text{old}} H_{\text{old}}^{-1} = I$, where $I$ is the identity matrix of appropriate dimensions, and therefore

$$\mathbf{b}_1^T \mathbf{h}_1 + bh + \mathbf{b}_2^T \mathbf{h}_2 = 1 > 0, \qquad (46)$$

where the dot product of the $i$-th row of $H_{\text{old}}$ in (44) and the $i$-th column of $H_{\text{old}}^{-1}$ in (37a) is taken and 1 is the $(i,i)$-th element of $I$. By neglecting the 1 in (46), multiplying both sides by $\lambda = \lambda_{\min} < 0$ from (27) or (28) (which flips the inequality sign) and dividing both sides by $h > 0$ (which is the $(i,i)$-th element in $H_{\text{old}}^{-1}$ in (37a) and $H_{\text{old}} \succ 0$ from Assumption 1.2), we can rewrite this inequality as

$$b\lambda < -\frac{\lambda}{h} \mathbf{b}_1^T \mathbf{h}_1 - \frac{\lambda}{h} \mathbf{b}_2^T \mathbf{h}_2. \qquad (47)$$

Now, $\mathbb{A}_{\text{new}}^m = \mathbb{A}_{\text{old}}^m \setminus \{\mathbb{A}_i\}$ and the $\mathbb{A}_i$-th row in (29) gives

$$K_{\mathbb{A}_i} = k + \mathbf{b}_1^T \boldsymbol{\lambda}_1 + \mathbf{b}_2^T \boldsymbol{\lambda}_2 - \frac{\lambda}{h} \mathbf{b}_1^T \mathbf{h}_1 - \frac{\lambda}{h} \mathbf{b}_2^T \mathbf{h}_2 > 0, \quad (48)$$

where the dot product of the $i$-th row (excluding the $i$-th column) of $H_{\text{old}}$ in (44) and $\lambda_{\text{new}}$ in (39) is taken. The inequality in (48) follows from substituting $b\lambda$ in (47) into (45). Since $K_{\mathbb{A}_i} \not< 0$, the **while** condition in step 4 of Algorithm 1 ensures constraint index $\mathbb{A}_i$ is not added again to $\mathbb{A}_{\text{old}}^{m+1} = \mathbb{A}_{\text{new}}^m$ at iteration $m+1$. ∎

We arrive at the following theorem.

**Theorem 1.** *Suppose Assumption 1.1 and Assumption 1.2 hold. Then the imuQP algorithm in Algorithm 1 is guaranteed to terminate within a finite number of iterations.*

*Proof.* From Lemma 1.1, Lemma 1.2.(ii), Lemma 1.3 and the trivial fact that after removing constraint index $\mathbb{A}_i$ from $\mathbb{A}_{\text{old}}^m$ such that $\mathbb{A}_{\text{new}}^m = \mathbb{A}_{\text{old}}^m \setminus \{\mathbb{A}_i\}$, it cannot be removed again from $\mathbb{A}_{\text{old}}^{m+1} = \mathbb{A}_{\text{new}}^m$, infinite cycling is prevented (i.e. repeatedly adding and removing the same constraints). Moreover, since the number of constraints $p$ in (13) is finite, and adding and removing constraints from the active set will result in $K$ in (29) and $\lambda_{\mathbb{A}}$ consist of non-negative numbers, the **while** loops in steps 4 and 4.d in Algorithm 1, respectively, will be exited within a finite number of iterations. ∎

*Remark 3.* When the unconstrained solution $\theta_0$ violates no or one constraint, imuQP terminates after one or two iterations and $c^*$ is equal to zero or one, respectively, because imuQP adds the most violated constraint index to the active set.

*2) Feasibility and optimality of the constrained solution:* Next, a set of instrumental lemmas, which makes up the proof of the main result – feasibility (related to the feasible set $\mathcal{X}$) and optimality (related to $\bar{J}(\theta)$ in (11)) of the constrained solution, is presented. The following lemma considers the costs of the unconstrained and constrained solutions.

**Lemma 2.1.** *Given $\theta_0$ in (21) and $\theta^*$ in (31), then it holds that $\bar{J}(\theta_0) < \bar{J}(\theta^*)$ with $\bar{J}(\theta)$ in (11).*

*Proof.* Substituting $\theta_0$ and $\theta^*$ into $\bar{J}(\theta)$ in (11) results in

$$\bar{J}_0 = \bar{J}(\theta_0) = -\frac{1}{2} F^T E^{-1} F \tag{49a}$$

$$< \bar{J}_0 + \frac{1}{2} (\lambda_{\mathbb{A}}^*)^T H_{\mathbb{A}} \lambda_{\mathbb{A}}^* = \bar{J}(\theta^*), \tag{49b}$$

since $H_{\mathbb{A}} \succ 0$ (Assumption 1.2) and a non-zero vector $\lambda_{\mathbb{A}}^* \in \mathbb{R}_{>0}^c$ gives $(\lambda_{\mathbb{A}}^*)^T H_{\mathbb{A}} \lambda_{\mathbb{A}}^* > 0$ [42, proof obs. 7.1.9, p. 432]. ∎

*Remark 4.* The expression $\bar{J}(\theta_0) = -\frac{1}{2} F^T E^{-1} F$ in (49a) is also shown in [41, prop. 4.1, p. 9].

We proceed with two lemmas regarding which constraint index should be added or removed from the active set and how the cost function $\bar{J}(\theta)$ in (11) will change.

**Lemma 2.2.** *Suppose Assumption 1.2 holds and at iteration $m$ we have $\mathbb{A}_{\text{new}}^m$ and the solution $\theta_m$ violates the constraint with index $j$ the most and corresponds to $K_{\min} < 0$. Then at iteration $m+1$, (i) constraint index $j$ should be added to $\mathbb{A}_{\text{old}}^{m+1} = \mathbb{A}_{\text{new}}^m$ such that $\mathbb{A}_{\text{new}}^{m+1} = \mathbb{A}_{\text{old}}^{m+1} \cup \{j\}$ (step 4.a.ii in Algorithm 1), and (ii) for the solution $\theta_{m+1}$ it holds that $\bar{J}(\theta_{m+1}) > \bar{J}(\theta_m)$.*

*Proof.* (i) $K_{\min} < 0$ is the $j$-th and most negative element of $K_0$ in (22) (or $K$ in (29)) and equals $\gamma_j - m_j \theta_0$ (or $\gamma_j - m_j \theta_m$). It indicates which constraint $m_j \theta_0 \leq \gamma_j$ (or $m_j \theta_m \leq \gamma_j$) is violated the most by the solution $\theta_0$ (or $\theta_m$) with $\mathbb{A}_{\text{new}}^m$. By adding constraint index $j$ to $\mathbb{A}_{\text{old}}^{m+1}$ such that $\mathbb{A}_{\text{new}}^{m+1} = \mathbb{A}_{\text{old}}^{m+1} \cup \{j\}$, most of the other violated constraints are also covered and it holds that $m_j \theta_0 = \gamma_j$ (or $m_j \theta_m = \gamma_j$) (Lemma 1.2.(i)).

(ii) Case 1 (from Lemma 1.1.(i)): substituting $H_{\mathbb{A}} = h_{jj}$ and $\lambda_{\mathbb{A}}$ in (24) into $\bar{J}(\theta^*)$ in (49b), with $\bar{J}_0$ in (49a), gives

$$\bar{J}(\theta_1) = \bar{J}_0 + \frac{1}{2} \frac{K_{\min}^2}{h_{jj}} > -\frac{1}{2} F^T E^{-1} F = \bar{J}(\theta_0) = \bar{J}_0,$$

where $h_{jj} > 0$ (because $H_{\mathbb{A}} \succ 0$, see Assumption 1.2) and $K_{\min} < 0$ in (23). This result also follows from Lemma 2.1.

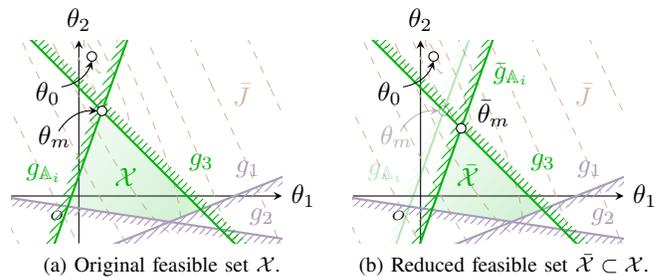

(a) Original feasible set $\mathcal{X}$.    (b) Reduced feasible set $\bar{\mathcal{X}} \subset \mathcal{X}$.

Fig. 1. Cost function $\bar{J}$, active constraints $g_{\mathbb{A}_i}$ and $g_3$ and inactive constraints $g_1$ and $g_2$, where $g_{\mathbb{A}_i}$ is perturbed inwards the interior of $\mathcal{X}$.

Case 2 (from Lemma 1.1.(i)): substituting $H_{\text{new}}$ in (33) and $\lambda_{\text{new}}$ in (36) into $\bar{J}(\theta^*)$ in (49b), with $\bar{J}_0$ in (49a), gives

$$\bar{J}(\theta_{m+1}) = \bar{J}_0 + \frac{1}{2} \lambda_{\text{new}}^T H_{\text{new}} \lambda_{\text{new}}$$

$$= \bar{J}_0 + \frac{1}{2} \lambda_{\text{old}}^T H_{\text{old}} \lambda_{\text{old}} + \frac{1}{2} \frac{\left(k + \mathbf{h}^T \lambda_{\text{old}}\right)^2}{h - \mathbf{h}^T H_{\text{old}}^{-1} \mathbf{h}}$$

$$> \bar{J}_0 + \frac{1}{2} \lambda_{\text{old}}^T H_{\text{old}} \lambda_{\text{old}} = \bar{J}(\theta_m).$$

From Assumption 1.2 the $j$-th constraint is linearly independent on the $n$-th constraints with $n \in \mathbb{A}_{\text{old}}^{m+1}$. Therefore $h - \mathbf{h}^T H_{\text{old}}^{-1} \mathbf{h} > 0$, as will be shown in Lemma 3.1.(ii). From iteration $m$, we have that $K_{\min} = k + \mathbf{h}^T \lambda_{\text{old}} < 0$ in (30). ∎

**Lemma 2.3.** *Suppose at iteration $m$ we have $\mathbb{A}_{\text{new}}^m$ and the solution $\theta_m$ causes the constraint with index $\mathbb{A}_i$, the $i$-th entry in $\mathbb{A}_{\text{new}}^m$, to become the most inactive and corresponds to $\lambda_{\min} < 0$. Then at iteration $m+1$, (i) constraint index $\mathbb{A}_i$ should be removed from $\mathbb{A}_{\text{old}}^{m+1} = \mathbb{A}_{\text{new}}^m$ such that $\mathbb{A}_{\text{new}}^{m+1} = \mathbb{A}_{\text{old}}^{m+1} \setminus \{\mathbb{A}_i\}$ (step 4.d.ii in Algorithm 1), and (ii) for the solution $\theta_{m+1}$ it holds that $\bar{J}(\theta_{m+1}) < \bar{J}(\theta_m)$.*

*Proof.* (i) Suppose the opposite is true and $\theta_m$ is the optimal and feasible solution to the QP in (11) with the Lagrange multiplier $\lambda$. We have that $\lambda_{\min} = \lambda_{\mathbb{A}_i} < 0$, which is the $i$-th and most negative element from $\lambda_{\mathbb{A}}$ in (27) or (28). From Lemma 1.2.(i), it holds that $m_{\mathbb{A}_i} \theta_m = \gamma_{\mathbb{A}_i}$, see Fig. 1a. Substituting (19) into (20) with $\theta_m$ instead of $\theta_0$ gives

$$\nabla \bar{J}(\theta_m) = -M^T \lambda.$$

For small $\delta \bar{J}(\theta_m)$ and $\delta \theta$ it holds that $\nabla \bar{J}(\theta_m) \approx \frac{\delta \bar{J}(\theta_m)}{\delta \theta}$ and

$$\delta \bar{J}(\theta_m) = \delta \theta \cdot \nabla \bar{J}(\theta_m) = -\delta \gamma^T \lambda = -\sum_{n \in \mathbb{A}_{\text{old}}^{m+1}} \delta \gamma_n \lambda_n,$$

where $M \delta \theta = \delta \gamma$ and $\lambda_n = 0$ for $n \notin \mathbb{A}_{\text{old}}^{m+1}$ from (14b). Let us perturb the active constraint $g_{\mathbb{A}_i}$ inwards the interior of $\mathcal{X}$. Let $\bar{\theta}_m = \theta_m + \delta \theta$ be the solution to the modified QP which still satisfies the perturbed active constraint $\bar{g}_{\mathbb{A}_i}$ with equality

$$\bar{g}_{\mathbb{A}_i}: m_{\mathbb{A}_i} \bar{\theta}_m = m_{\mathbb{A}_i} (\theta_m + \delta \theta) = \gamma_{\mathbb{A}_i} + \delta \gamma_{\mathbb{A}_i} < \gamma_{\mathbb{A}_i},$$

where $m_{\mathbb{A}_i} \delta \theta = \delta \gamma_{\mathbb{A}_i} < 0$. Therefore, $\bar{\theta}_m \in \bar{\mathcal{X}} \subset \mathcal{X}$, where $\bar{\mathcal{X}}$ is the reduced feasible set, see Fig. 1b. Since the other active constraints remain unchanged, i.e. $\delta \gamma_n = 0$ for $n \in \mathbb{A}_{\text{old}}^{m+1} \setminus \{\mathbb{A}_i\}$, we have $\delta \bar{J}(\theta_m) = -\delta \gamma_{\mathbb{A}_i} \lambda_{\mathbb{A}_i} < 0$ and $\bar{J}(\bar{\theta}_m) < \bar{J}(\theta_m)$. Therefore, $\theta_m$ is not a minimizer to the original QP in (11), which is a contradiction. By removing constraint index $\mathbb{A}_i$ from $\mathbb{A}_{\text{old}}^{m+1}$ such that $\mathbb{A}_{\text{new}}^{m+1} = \mathbb{A}_{\text{old}}^{m+1} \backslash \{\mathbb{A}_i\}$, $\bar{J}(\theta)$ can be decreased the most while staying inside $\mathcal{X}$ [35, ch. 12.3, p. 362-363] [43, ch. 5.6, p. 249-253] and we get much closer to the minimizer. The elements of the Lagrange multiplier $\lambda$ are also known as *shadow prices* [43, ch. 5.4.4, p. 241].

*(ii)* From Lemma 1.1.(i), constraints cannot be removed from $\mathbb{A}_{\text{old}}^{m+1}$ at iterations $m = 1, 2, 3$. Therefore, we consider iteration $m \geq 4$ with $|\mathbb{A}_{\text{old}}^m| > 1$, and from (35) we have that $H_{\text{old}}\lambda_{\text{old}} = -\mathbf{k}_{\text{old}}$ and $H_{\text{new}}\lambda_{\text{new}} = -\mathbf{k}_{\text{new}}$. Using $\bar{J}(\theta^*)$ in (49b), $\lambda_{\text{old}}$ and $\mathbf{k}_{\text{old}}$ in (37b) and $\lambda_{\text{new}}$ and $\mathbf{k}_{\text{new}}$ in (39) gives

$$\begin{aligned}\bar{J}(\theta_{m+1}) &= \bar{J}_0 + \frac{1}{2}\lambda_{\text{new}}^T H_{\text{new}}\lambda_{\text{new}} = \bar{J}_0 - \frac{1}{2}\lambda_{\text{new}}^T \mathbf{k}_{\text{new}} \\ &= \bar{J}_0 + \frac{1}{2}\left(\frac{\lambda}{h}\mathbf{h}_1^T\mathbf{k}_1 + \frac{\lambda}{h}\mathbf{h}_2^T\mathbf{k}_2 - \boldsymbol{\lambda}_1^T\mathbf{k}_1 - \boldsymbol{\lambda}_2^T\mathbf{k}_2\right) \\ &< \bar{J}_0 + \frac{1}{2}\left(-\lambda k - \boldsymbol{\lambda}_1^T\mathbf{k}_1 - \boldsymbol{\lambda}_2^T\mathbf{k}_2\right) \\ &= \bar{J}_0 - \frac{1}{2}\lambda_{\text{old}}^T\mathbf{k}_{\text{old}} = \bar{J}_0 + \frac{1}{2}\lambda_{\text{old}}^T H_{\text{old}}\lambda_{\text{old}} = \bar{J}(\theta_m),\end{aligned}$$

where we used $\bar{J}_0$ in (49a). To prove this inequality, note that the $i$-th row of $\lambda_{\text{old}} = -H_{\text{old}}^{-1}\mathbf{k}_{\text{old}}$ in (35) gives

$$\lambda = -\mathbf{h}_1^T\mathbf{k}_1 - hk - \mathbf{h}_2^T\mathbf{k}_2 < 0, \quad (50)$$

where the dot product of the $i$-th row of $H_{\text{old}}^{-1}$ in (37a) and $\mathbf{k}_{\text{old}}$ in (37b) is taken and $\lambda = \lambda_{\min} < 0$ is the $i$-th element of $\lambda_{\text{old}}$ in (37b). By neglecting $\lambda$ in (50), multiplying both sides by $\lambda = \lambda_{\min} < 0$ from (27) or (28) (which flips the inequality sign) and dividing both sides by $h > 0$ (which is the $(i, i)$-th element in $H_{\text{old}}^{-1}$ in (37a) and $H_{\text{old}} \succ 0$ from Assumption 1.2), we can rewrite this inequality as

$$\frac{\lambda}{h}\mathbf{h}_1^T\mathbf{k}_1 + \frac{\lambda}{h}\mathbf{h}_2^T\mathbf{k}_2 < -\lambda k. \quad \blacksquare$$

We can conclude with the following theorem.

**Theorem 2.** *Suppose Assumption 1.1 and Assumption 1.2 hold. Then the solution computed by the imuQP algorithm in Algorithm 1 is feasible and optimal.*

*Proof.* As was shown in Theorem 1, the imuQP algorithm in Algorithm 1 terminates within a finite number of iterations. Let $\theta^*$ denote the solution obtained after the last iteration. Adding constraints to the active set ensures $\theta^*$ is feasible (Lemma 1.2 and Lemma 2.2). Removing constraints from the active set ensures $\theta^*$ is both feasible (Lemma 1.3) and optimal (Lemma 2.3). The two **while** conditions in steps 4 and 4.d in Algorithm 1 ensure $\theta^*$ is both feasible and optimal using $K = \gamma - M\theta^* \geq 0$ and $\lambda^* \geq 0$ in (15), respectively. $\blacksquare$

*Remark 5.* In imuQP, $\theta_0$ in (19) and $\theta^*$ in (17) are used in steps 1 and 5 in Algorithm 1, respectively. AIMU and SIMU efficiently compute $\lambda_{\mathbb{A}}$ and $H_{\mathbb{A}}^{-1}$ based on (17). From (20), the constrained solution is the point where $\mathcal{X}$ touches a level curve of $\bar{J}(\theta)$, i.e. their gradients are parallel to each other.

*3) Linear dependence of constraints:* Finally, the following lemma shows how to check whether the to-be-added constraint is linearly dependent on the constraints present in the current active set — which is part of the inverse matrix update process.

**Lemma 3.1.** *Suppose at iteration $m$ we have $q$ in (25) (step 4.b.i in Algorithm 1) and $n \in \mathbb{A}_{\text{old}}^m$, which is the current active set. Then the to-be-added $j$-th constraint is either (i) linearly dependent on the $n$-th constraints if $q = 0$, or (ii) linearly independent on the $n$-th constraints if $q > 0$.*

*Proof.* (i) Suppose constraint index $j$ is added to $\mathbb{A}_{\text{old}}^m$ such that $\mathbb{A}_{\text{new}}^m = \mathbb{A}_{\text{old}}^m \cup \{j\}$, $M_{\text{old}}$ (or $M_{\text{new}}$) consists of the $n$-th rows of $M$ with $n \in \mathbb{A}_{\text{old}}^m$ (or $n \in \mathbb{A}_{\text{new}}^m$) and $M_{\text{old}}$ has full row rank. Suppose $m_j$ can be written as a weighted sum of the rows of $M_{\text{old}}$, with weights $y$, i.e.

$$m_j = y^T M_{\text{old}} \quad \text{and} \quad M_{\text{new}} = \begin{bmatrix} M_{\text{old}} \\ m_j \end{bmatrix}. \quad (51)$$

From $H_{\mathbb{A}}$ in (18) and $H_{\text{new}}$ in (33), we can write

$$\begin{aligned}H_{\text{new}} &= M_{\text{new}}E^{-1}M_{\text{new}}^T \quad &(52)\\ &= \begin{bmatrix} M_{\text{old}}E^{-1}M_{\text{old}}^T & M_{\text{old}}E^{-1}m_j^T \\ m_jE^{-1}M_{\text{old}}^T & m_jE^{-1}m_j^T \end{bmatrix} = \begin{bmatrix} H_{\text{old}} & \mathbf{h} \\ \mathbf{h}^T & h \end{bmatrix}.\end{aligned}$$

From (52) and $m_j$ in (51), we can show that the weights $y$ in (25) can be computed as

$$H_{\text{old}}^{-1}\mathbf{h} = \left(M_{\text{old}}E^{-1}M_{\text{old}}^T\right)^{-1}\left(M_{\text{old}}E^{-1}M_{\text{old}}^T\right)y = y,$$

and using $m_j$ in (51) we can indeed show that $q$ in (25)

$$q = h - \mathbf{h}^T y = y^T\left(M_{\text{old}}E^{-1}M_{\text{old}}^T - M_{\text{old}}E^{-1}M_{\text{old}}^T\right)y = 0.$$

*(ii)* Being linearly independent means that $m_j$ cannot be written as a weighted sum of the rows of $M_{\text{old}}$ and both $M_{\text{new}}$ in (51) and $M_{\text{old}}$ have full row rank. Since $E \succ 0$, then $H_{\text{new}} = M_{\text{new}}E^{-1}M_{\text{new}}^T \succ 0$ in (52) and $H_{\text{old}} = M_{\text{old}}E^{-1}M_{\text{old}}^T \succ 0$ [42, obs. 7.1.8(b), p. 431] and from Schur complement [41, prop. 2.2, p. 4] we have that $q$ in (25)

$$q = h - \mathbf{h}^T H_{\text{old}}^{-1}\mathbf{h} = m_j S m_j^T > 0,$$

where $S = E^{-1} - E^{-1}M_{\text{old}}^T H_{\text{old}}^{-1}M_{\text{old}}E^{-1}$. Since $m_j$ is non-zero, it also follows that $S \succ 0$ [42, eq. (7.1.1a), p. 429]. $\blacksquare$

*Remark 6.* In imuQP, the **if** condition in step 4.b.ii in Algorithm 1 is rather implemented as a tolerance of $q \leq 1 \cdot 10^{-11}$ for numerical precision reasons; setting this tolerance too loose (large) will decrease accuracy, while setting it too tight (small) will cause the solver to never terminate.

This leads us to the following theorem.

**Theorem 3.** *Suppose at iteration $m$ we have the to-be-added $j$-th constraint with $q = h - \mathbf{h}^T H_{\text{old}}^{-1}\mathbf{h} = 0$ and $y_{\max} = \max(y) > 0$ in (25) and the current active set $\mathbb{A}_{\text{old}}^m$. Then the constraint with index $\mathbb{A}_f$, the $f$-th entry in $\mathbb{A}_{\text{old}}^m$, should be removed from $\mathbb{A}_{\text{old}}^m$ before adding the $j$-th constraint to $\mathbb{A}_{\text{old}}^m$, where $f = \arg\max(y)$ (steps 4.b.ii.C until 4.b.ii.E in Algorithm 1).*

*Proof.* From Lemma 3.1.(i), $m_j$ is linearly dependent on the $n$-th rows of $M$ with $n \in \mathbb{A}_{\text{old}}^m$ and weights $y = H_{\text{old}}^{-1}\mathbf{h}$ in (25). Directly performing AIMU with constraint index $j$ is not possible as this requires division by 0 in computing $H_{\text{new}}^{-1}$ in (34) and $\lambda_{\text{new}}$ in (36). From Lemma 2.3.(i), we look for the most negative element in $\lambda_{\text{old}}$ in order to remove its corresponding constraint index from $\mathbb{A}_{\text{old}}^m$. From Lemma 1.1, we should only consider the top expression of $\lambda_{\text{new}}$ in (36). Since the second term dominates and $K_{\min} = k + \mathbf{h}^T \lambda_{\text{old}} < 0$ in (29) from iteration $m-1$, this is equivalent to finding $y_{\max} = \max(H_{\text{old}}^{-1}\mathbf{h}) = \max(y) > 0$ with $f = \arg\max(y)$ in (26) and $y$ should have at least one positive element. ∎

*Remark 7.* For implementation reasons in Algorithm 1, constraint index $j$ is added to the list of active constraints in step 4.a.ii before removing constraint index $\mathbb{A}_f$ in step 4.b.ii.C. However, SIMU with constraint index $\mathbb{A}_f$ removed in step 4.b.ii.D are performed before AIMU with constraint index $j$ added in step 4.b.ii.E.

### A. Infeasibility detection of a QP

So far, we assumed that the QP in (11) is feasible. When a QP is infeasible, $\mathcal{X} = \varnothing$ and there is no solution that satisfies all inequality constraints. In literature, infeasibility of a QP with affine inequality constraints is usually checked using Farkas' lemma [44, prop. 6.4.3 and remarks, p. 90-92]. In imuQP, infeasiblity of a QP can easily be detected *for free* as it is part of the inverse matrix update process using $y$ in (25), as will be shown in the following theorem.

**Theorem 4.** *Suppose at iteration $m$ we have the to-be-added $j$-th constraint with $q = 0$ and $y_{\max} = \max(y) \leq 0$ in (25) and the current active set $\mathbb{A}_{\text{old}}^m$. Then the QP in (11) is infeasible (step 4.b.ii.B in Algorithm 1).*

*Proof.* Consider the following affine inequality constraints

$$g_1: \quad \theta_2 \leq a_1(\theta_1 - \theta_1^*) + \theta_2^* \tag{53a}$$
$$g_2: \quad \theta_2 \leq -a_2(\theta_1 - \theta_1^*) + \theta_2^* \tag{53b}$$
$$g_3: \quad \theta_2 \geq -\frac{-a_1 y_1 + a_2 y_2}{y_1 + y_2}(\theta_1 - \theta_1^*) + \theta_2^* + \epsilon, \tag{53c}$$

where $a_1, a_2, \epsilon > 0$, which can be written in the matrix form

$$\begin{bmatrix} m_1 \\ m_2 \\ m_3 \end{bmatrix} \theta \leq \begin{bmatrix} \gamma_1 \\ \gamma_2 \\ \gamma_3 \end{bmatrix}, \quad \text{where } \theta = \begin{bmatrix} \theta_1 \\ \theta_2 \end{bmatrix},$$

with

$$\begin{bmatrix} m_1 \\ m_2 \end{bmatrix} = \begin{bmatrix} -a_1 & 1 \\ a_2 & 1 \end{bmatrix}, \begin{bmatrix} \gamma_1 \\ \gamma_2 \end{bmatrix} = \begin{bmatrix} -a_1\theta_1^* + \theta_2^* \\ a_2\theta_1^* + \theta_2^* \end{bmatrix},$$
$$m_3 = \begin{bmatrix} -a_1 y_1 + a_2 y_2 & y_1 + y_2 \end{bmatrix}, \tag{54}$$
$$\gamma_3 = (-a_1 y_1 + a_2 y_2)\theta_1^* + (y_1 + y_2)(\theta_2^* + \epsilon).$$

Suppose at iteration $m$ we have $\mathbb{A}_{\text{old}}^m = \{1, 2\}$ and the solution $\theta^* = (\theta_1^*, \theta_2^*)$. Constraints $g_i: m_i\theta \leq \gamma_i$ with $i = \{1, 2\}$ are constructed such that their borders go through the point $\theta^*$ as

$$m_i(\theta - \theta^*) \leq 0, \quad \text{with } \gamma_i = m_i\theta^*.$$

By construction, $\mathcal{X}_{1,2} = \{\theta \,|\, m_1\theta \leq \gamma_1 \wedge m_2\theta \leq \gamma_2\}$ is always below $\theta^*$, see Fig. 2. However, $\theta^*$ violates constraint $g_3$: $m_3\theta \leq \gamma_3$ the most and therefore $\mathbb{A}_{\text{new}}^m = \mathbb{A}_{\text{old}}^m \cup \{3\} = \{1, 2, 3\}$, from Lemma 2.2.(i), i.e.

$$m_3\theta^* > \gamma_3 \quad \text{and} \quad \begin{bmatrix} m_1 \\ m_2 \end{bmatrix} \theta^* = \begin{bmatrix} \gamma_1 \\ \gamma_2 \end{bmatrix}, \tag{55}$$

where the latter comes from Lemma 1.2.(i) and $m_3$ is a linear combination of $m_1$ and $m_2$ with weights $y_1, y_2 \leq 0$. Dividing by $y_1 + y_2 < 0$ flips the inequality sign in (53c) ($y_1$ and $y_2$ are never both 0 as $m_3$ and $\gamma_3$ in (54) would be zero), i.e.

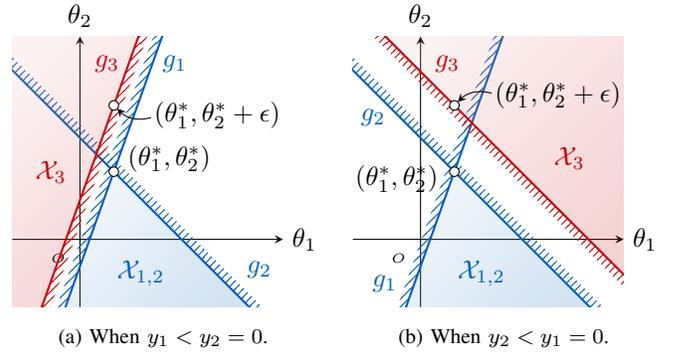

(a) When $y_1 < y_2 = 0$.  (b) When $y_2 < y_1 = 0$.

Fig. 2. Infeasibility of a QP where the violated constraint $g_3$ is linearly dependent on the active constraints $g_1, g_2$, with weights $y_1, y_2 \leq 0$.

$$m_3 = y^T \begin{bmatrix} m_1 \\ m_2 \end{bmatrix}, \quad \text{with } y = \begin{bmatrix} y_1 \\ y_2 \end{bmatrix}. \tag{56}$$

By substituting $m_3$ in (56) into (55), it holds that

$$m_3\theta^* = y^T \begin{bmatrix} m_1 \\ m_2 \end{bmatrix} \theta^* = y^T \begin{bmatrix} \gamma_1 \\ \gamma_2 \end{bmatrix} > \gamma_3,$$

and constraint $g_3$: $m_3\theta \leq \gamma_3$ is constructed as

$$m_3\left(\theta - \begin{bmatrix} \theta_1^* \\ \theta_2^* + \epsilon \end{bmatrix}\right) \leq 0, \quad \text{with } \gamma_3 = m_3 \begin{bmatrix} \theta_1^* \\ \theta_2^* + \epsilon \end{bmatrix}.$$

By construction, the border of constraint $g_3$ goes through the point $(\theta_1^*, \theta_2^* + \epsilon)$ and $\mathcal{X}_3 = \{\theta \,|\, m_3\theta \leq \gamma_3\}$ is always above this point. Constraint $g_3$ is therefore always violated by $\theta^*$, see Fig. 2 for some special cases of $y$. Clearly, $\mathcal{X} = \mathcal{X}_{1,2} \cap \mathcal{X}_3 = \varnothing$ for $y_1, y_2 \leq 0$ and therefore the QP is infeasible. It can be generalized that if $y_{\max} = \max(y) \leq 0$, i.e. all elements in $y$ in (25) are non-positive, then the QP in (11) is infeasible. ∎

### B. Memory footprint

The memory footprint gives a rough measure of how many numbers need to be stored, where minor contributions of scalar numbers are ignored. To simplify the analysis in what follows, we take $n_u = n_y = 1$ and $p$ and $w$ in (13) become

$$p = 6N + 2 \quad \text{and} \quad w = \tfrac{p}{2} = 3N + 1. \tag{57}$$

*1) imuQP:* The Hessian matrix $E$ is manually factorized offline using Gaussian elimination (without pivoting) [45, ch. 3.2.6, p. 98] into strictly lower and not strictly upper triangular matrices $E_L$ and $E_U$, respectively, which is equivalent to

$$E_L = L - I, \quad E_U = U, \quad E = LU = (E_L + I)E_U,$$

where $L$ and $U$ are lower and upper triangular matrices from the LU decomposition – obtained from the MATLAB function lu() – and $I$ is the identity matrix of appropriate dimensions. It is worth mentioning that generating $E_L$ and $E_U$ using manual Gaussian elimination or lu() gives a negligible difference in the solution of imuQP, although the former method is slightly more accurate and is therefore used in imuQP. $E_L$ and $E_U$ are used to perform forward and backward substitution, respectively [45, ch. 3.1.1-3.1.2, p. 88-89] to compute offline $H$ in (32), online $\theta_0$ in (21) and online $\theta^*$ in (31). Factorization is done once manually offline, since $E$ stays the same while only the right-hand side changes. Since $E_L$ and $E_U$ do not overlap, they can be safely added up and stored in a single dense matrix $E_{LU} = E_L + E_U$ [45, ch. 3.2.8, p. 99]. Even though Algorithm 1 takes as input the Hessian matrix $E$, in the actual implementation $E_{LU}$ is given as input.

*2) qpOASES:* From [8, ch. 4.1, p. 39-42], [8, ch. 4.6.2 and tab. 4.4, p. 66] and [7, ch. 4.1.4, p. 54], both $w$ and $\Delta w$ have the dimensions of $x(k) \in \mathbb{R}^{n_x}$, $\mathcal{R}_k \in \mathbb{R}^{N n_y}$ and $u(k-1) \in \mathbb{R}^{n_u}$ added up as they appear as parameters in $F$ and $\gamma$ in (11). Using $n_u = n_x = n_y = 1$, we have that both $w, \Delta w \in \mathbb{R}^{N+2}$. We ignore the terms with $\mathcal{O}(n)$, where $n$ is the number of decision variables in qpOASES.

*3) Comparison:* Table II summarizes how many numbers need to be stored by imuQP and qpOASES. It can be concluded that imuQP (with leading term $52N^2$) has to store more numbers in total than qpOASES (with leading term $10N^2$), but both are within the same order of complexity, i.e. $\mathcal{O}\left(N^2\right)$.

### C. Computational complexity

The computational complexity gives a rough measure of how many arithmetic operations are performed online such as additions, subtractions, multiplications and divisions, where minor contributions are ignored. We again take $n_u = n_y = 1$, where $p$ is defined in (57). Computing $AB = C \in \mathbb{R}^{m \times p}$ with $A \in \mathbb{R}^{m \times n}$ and $B \in \mathbb{R}^{n \times p}$ takes $mp(2n-1)$ operations, as each element in $C$ takes $n$ multiplications and $n-1$ additions.

*1) imuQP:* Offline manual factorization using Gaussian elimination (without pivoting) takes $\frac{2}{3}N^3$ operations [45, ch. 3.2.6, p. 98]. Forward and backward substitution each takes $N^2$ operations and thus $2N^2$ operations in total [45, ch. 3.1.1-3.1.2, p. 88-89], to compute offline $H$ in (32), online $\theta_0$ in (21) and online $\theta^*$ in (31). Since $E$ stays the same while only the right-hand side changes, expensive factorization only needs to be done once offline, while cheap forward and backward substitutions can be performed online. Furthermore, since $M$ in (11) stays the same, $H$ in (32) is computed offline, which takes $2pN^2 + p^2(2N-1)$ operations. Table III gives an analysis of how many arithmetic operations are performed online by imuQP and the steps refer to the ones in Algorithm 1. Here, $t_l \in \mathbb{Z}_{\geq 0}$ ($t_a \in \mathbb{Z}_{\geq 0}$) denotes how many times a constraint is added to the current active set $\mathbb{A}$ which is linearly dependent (independent) on the constraints in $\mathbb{A}$, respectively, $t_r \in \mathbb{Z}_{\geq 0}$ denotes how many times a constraint is removed from $\mathbb{A}$. The total number of arithmetic operations $\mathbf{n}$ is then

$$\mathbf{n} = 4N^2 + N\left[2(p+c)+1\right] + p + 2 + t_l\left[10c^2 + 20\tfrac{1}{2}c + 2p(c+1) + 6\right] + t_a\left[4c^2 + 9c + 2p(c+1) + 3\right] + t_r\left[4c^2 + 6\tfrac{1}{2}c + 1\right].$$

Since usually $1 \ll c \ll p$, the following terms are dominating

$$\bar{\mathbf{n}} = 4N^2 + N(2p+1) + p + t_l\left[10c^2 + 2p(c+1)\right] + t_a\left[4c^2 + 2p(c+1)\right] + t_r\left[4c^2\right].$$

*2) qpOASES:* From [8, ch. 4.6.1, p. 62-63] and [7, ch. 4.1.4, p. 53-54], the number of active constraints $n_\mathbb{A}$, the number of inactive constraints $n_\mathbb{I}$ and the number of decision variables $n$ from qpOASES correspond to $c$, $p-c$ and $N$ from imuQP, respectively. Again, we ignore the terms with $\mathcal{O}(n)$ and $n_\mathbb{X}$ as we only consider operations involving adding and removing affine inequality constraints from the active set $\mathbb{A}$ and ignore those that treat box constraints separately. Then $\mathbf{n}$ is

$$\mathbf{n} = t_l\left[10N^2 + 2\tfrac{1}{2}c^2 + N(p-6c)\right] + t_a\left[10N^2 + 2c^2 + N(p-7c)\right] + t_r\left[7\tfrac{1}{2}N^2 + 2\tfrac{7}{8}c^2 + N(p-3\tfrac{1}{2}c)\right],$$

which cannot be simplified, even if we use that $1 \ll c \ll p$.

TABLE II
MEMORY FOOTPRINT OF IMUQP (TOP) AND QPOASES (BOTTOM).

|  |  | Input data | Internal data |
|---|---|---|---|
| imuQP | Data | $E_{LU} \in \mathbb{R}^{N\times N}, F \in \mathbb{R}^N,$ $M \in \mathbb{R}^{p\times N}, \gamma \in \mathbb{R}^p,$ $H \in \mathbb{R}^{p\times p}$ | $\theta^* \in \mathbb{R}^N,$ $K, K_0 \in \mathbb{R}^p,$ $H_\mathbb{A}^{-1} \in \mathbb{R}^{w\times w},$ $\lambda_\mathbb{A}, \mathbb{A} \in \mathbb{R}^w$ |
|  | General | $N^2 + N + pN + p + p^2$ | $N + 2p + w^2 + 2w$ |
|  | In terms of $N$ | $43N^2 + 33N + 6$ | $9N^2 + 25N + 7$ |
|  | Total | $52N^2 + 58N + 13$ | |
| qpOASES | Data | $E \in \mathbb{R}^{N\times N}, F \in \mathbb{R}^N,$ $M \in \mathbb{R}^{p\times N}, \gamma \in \mathbb{R}^p$ | $w, \Delta w \in \mathbb{R}^{N+2},$ $\Delta F, \theta, \Delta \theta \in \mathbb{R}^N,$ $\Delta\gamma, \lambda, \Delta\lambda, \mathbb{A} \in \mathbb{R}^p,$ $R, T, Q \in \mathbb{R}^{N\times N}$ |
|  | General | $N^2 + N + pN + p$ | $3N^2 + 5N + 4p + 4$ |
|  | In terms of $N$ | $7N^2 + 9N + 2$ | $3N^2 + 29N + 12$ |
|  | Total | $10N^2 + 38N + 14$ | |

*3) Comparison:* Comparing $\mathbf{n}$ for both QP solvers, it can be concluded that imuQP (with leading terms $4N^2$ and $(10t_l + 4t_a + 4t_r)c^2$) has to perform fewer arithmetic operations than qpOASES (with leading terms $(10t_l + 10t_a + 7\tfrac{1}{2}t_r)N^2$ and $(2\tfrac{1}{2}t_l + 2t_a + 2\tfrac{7}{8}t_r)c^2$), but again with the same order of complexity, i.e. $\mathcal{O}\left(N^2 + c^2\right)$.

## V. SIMULATIONS: BENCHMARK EXAMPLE

SPEED and accuracy of imuQP will be compared with state-of-the-art active-set QP solvers by means of a MAT-LAB simulation. We will use a chain of six masses from [46], which is a well-known benchmark example and its diagram is shown in Fig. 3. A force $\mathcal{F}_i$ is acting on each mass $m_i$, with $i = 1, \ldots, 6$, which is connected to two springs with spring constants $k_i$ and $k_{i+1}$. Here, $\rho_i$ denotes the horizontal position of mass $m_i$ with respect to its rest position and for the two solid walls we have $\rho_0 = \rho_7 = 0$. Newton's second law gives

$$m_i \ddot{\rho}_i = -k_i(\rho_i - \rho_{i-1}) + k_{i+1}(\rho_{i+1} - \rho_i) + \mathcal{F}_i,$$

TABLE III
COMPUTATIONAL COMPLEXITY OF IMUQP. THE ABBREVIATIONS ARE: OCC.: NUMBER OF OCCURRENCES; NR.: NUMBER; OPERS.: OPERATIONS.

| Occ. | Step | Nr. of opers. | Occ. | Step | Nr. of opers. |
|---|---|---|---|---|---|
| 1 | 1 | $2N^2$ | $t_a$ | 4.b.ii | $2c^2 + 5c + 2$ |
|  | 2 | $2pN$ | $t_l + t_a$ | 4.c | $c$ |
|  | 3.b | $p$ |  | 4.d.ii | $\frac{c}{2}$ |
|  | 4.b | $2$ | $t_r$ | 4.d.iii | $4c^2 + 5c + 1$ |
| $t_l + t_a$ | 4.b.i | $2c^2 + 3c + 1$ |  | 4.d.iv | $c$ |
| $t_l$ | 4.b.ii.A | $c$ | $t_l + t_a$ | 4.e | $p(2c+1)$ |
|  | 4.b.ii.C | $\frac{c}{2}$ |  | 4.f | $p$ |
|  | 4.b.ii.D | $4c^2 + 5c + 1$ | 1 | 5 | $2N^2 +$ $N(2c+1)$ |
|  | 4.b.ii.E | $4c^2 + 10c + 4$ |  |  |  |

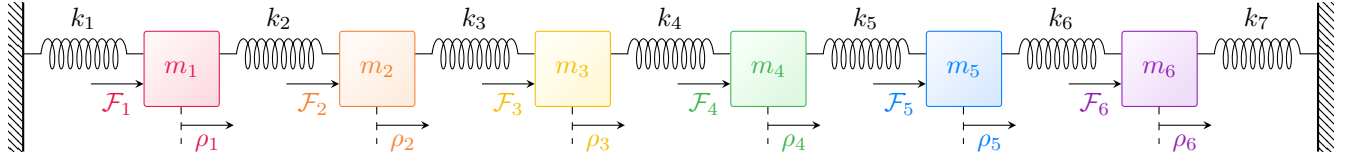

Fig. 3. Six masses $m_i$ connected to springs with constants $k_i$ between two solid walls, with forces $\mathcal{F}_i$, horizontal positions $\rho_i$ and $i = 1, \ldots, 6$.

and for simplicity, setting $m_i = k_i = 1$ gives

$$\ddot{\rho}_i = \rho_{i-1} - 2\rho_i + \rho_{i+1} + \mathcal{F}_i. \quad (58)$$

Now, with $n_u = n_y = 6$, we rewrite (58) as the continuous-time state-space model $\Sigma_c$ in (1) with $x_p = \begin{bmatrix} \rho^T & \dot{\rho}^T \end{bmatrix}^T$, $u = \mathcal{F}$, $y = \rho$, $w_c = 0$, but where '(t)' has been omitted, and

$$A_c = \begin{bmatrix} 0 & I \\ \mathcal{M} & 0 \end{bmatrix}, B_c = \begin{bmatrix} 0 \\ I \end{bmatrix}, C_c = \begin{bmatrix} I & 0 \end{bmatrix}, D_c = 0,$$

where $I$ is the $6 \times 6$ identity matrix and we define

$$\rho = \begin{bmatrix} \rho_1 \\ \vdots \\ \rho_6 \end{bmatrix}, \mathcal{F} = \begin{bmatrix} \mathcal{F}_1 \\ \vdots \\ \mathcal{F}_6 \end{bmatrix}, \mathcal{M} = \begin{bmatrix} -2 & 1 & 0 & 0 \\ 1 & \ddots & \ddots & 0 \\ 0 & \ddots & \ddots & 1 \\ 0 & 0 & 1 & -2 \end{bmatrix}.$$

From $\Sigma_c$ in (1) and the MATLAB function c2d() with $T_s = 4$ ms, the discrete-time state-space model $\Sigma_d$ in (2) is obtained. Here, $A_d$ is stable but not asymptotically stable (its discrete-time poles are located on the unit circle) since no dampers are present, the matrix pair $(A_d, B_d)$ is controllable and the matrix pair $(A_d, C_d)$ is observable. The augmented discrete-time state-space model for integral-action MPC in (3) is constructed and we impose the constraints in (8). Note that in [46], $\rho$ and $\dot{\rho}$ are limited to huge values of $\pm 1 \cdot 10^{15}$, while here we limit $y = \rho$ to $\pm 1 \cdot 10^6$ and impose no constraints on $\dot{\rho}$. In both cases $\rho$ and $\dot{\rho}$ are effectively unconstrained. The initial positions $\rho_i(0) \in [-1.2, 1.2]$ are randomly generated using rand() in MATLAB with a seed, initial velocities are $\dot{\rho}_i(0) = 0$ and each $\rho_i(t)$ tracks a sinusoidal reference $r_i(t)$

$$r_i(t) = A_r \sin\left(2\pi f_r t - 0.9 \cdot \frac{2\phi}{n_y - 1}\pi\right), \quad (59)$$

with $\phi = i - 1$. In what follows, $A_r$, $f_r$ and all MPC parameters are given in Table IV, where $I$ is the $6 \times 6$ identity matrix and $\bar{m}$ is the maximum number of iterations permitted in imuQP.

It is worth mentioning that ODYS and qpOASES call MEX functions. For a fair comparison, imuQP and Hildreth are written as MATLAB functions (.m files), from which MEX functions (.mexw64 files) and C code (.c files) are automatically generated using codegen from MATLAB Coder. Calling the MEX functions, which in turn call the C code, results in a much lower CPU time than directly calling the MATLAB functions. The simulations are run in MATLAB R2021b on an HP Envy x360 Convertible 13-bd0xxx device, which runs Microsoft Windows 11 Home 21H2, 11$^{\text{th}}$ Gen Intel® Core™ i7-1165G7 @2.80 GHz 2.80 GHz, 16.0 GB RAM (15.8 GB available) and 64-bit operating system, x64-based processor.

In the upcoming tables, the accuracy of each QP solver is assessed by evaluating the accuracy measures in Table V, at each time instant $k$. Here, $\theta^*$ and $\theta^*_{\text{quadprog}}$ are the solutions of some QP solver and quadprog, respectively, where the latter is considered optimal. The average value $\mathcal{T}_{\text{avg}}$ and maximum value $\mathcal{T}_{\max}$ are defined, respectively, as

$$\mathcal{T}_{\text{avg}} = \frac{1}{k_b - k_a + 1} \sum_{i=k_a}^{k_b} \tau_i, \quad \mathcal{T}_{\max} = \max_{i \in [k_a, k_b]} \tau_i, \quad (60)$$

where $k_a, k_b \in \mathbb{Z}_{>0}$ are the starting and final samples, respectively, with $k_a < k_b$ and $\tau_i$ is some data value at sample $i \in [k_a, k_b]$. The minimum value $\mathcal{T}_{\min}$ is defined similarly to $\mathcal{T}_{\max}$. For measuring CPU time in MATLAB, tic and toc() are called before and after the QP solver call, respectively.

TABLE IV
MPC PARAMETERS FOR SIMULATING A CHAIN OF SIX MASSES.

| Symbol | Value | Unit | Symbol | Value | Unit |
|---|---|---|---|---|---|
| $N$ | 27 | [–] | $u_{\min}$ | $-1$ | [N] |
| $T_s$ | 4 | [ms] | $u_{\max}$ | 1 | [N] |
| $Q$ | $210I$ | [–] | $y_{\min}$ | $-1 \cdot 10^6$ | [m] |
| $R$ | $0.008I$ | [–] | $y_{\max}$ | $1 \cdot 10^6$ | [m] |
| $P$ | $45Q$ | [–] | $\bar{m}$ | $\infty$ | [–] |
| $\Delta u_{\min}$ | $-0.5$ | [N] | $A_r$ | 1.2 | [m] |
| $\Delta u_{\max}$ | $0.5$ | [N] | $f_r$ | $\frac{1}{30}$ | [Hz] |

In what follows, the QP solver versions and settings used are summarized in Table VI. Here, lambda denotes $\lambda^*$, theta denotes $\theta^*$, m denotes the number of iterations performed $m^*$, act denotes the constraint indices of the optimal active set $\mathbb{A}^*$, c denotes $c^* = |\mathbb{A}^*|$ which is found by counting the number of non-zero elements in $\lambda^*$, lambda_act denotes $\lambda^*_{\mathbb{A}}$ and height_M denotes the number of rows in $M$. OSQP [2] was excluded in the comparison because it is not a pure active-set QP solver. However, it was shown in a recent paper [47] that OSQP's performance is comparable with qpOASES'.

### A. N = 27, single initial condition

A single simulation was run for 70 s − or 17500 samples with $T_s = 4$ ms − with a seed of 5 for the initial position $\rho_i(0)$. The trajectories of $y(t)$, $u(t)$ and $\Delta u(t)$, using imuQP as the QP solver, are shown in Fig. 4a-4c. Here, the reference signal $r(t)$ and the constraints for $u(t)$ and $\Delta u(t)$ are designed such that many constraints become active to push the QP solvers to the limit and it can be seen that these constraints are respected. Therefore, $y(t)$ may not reach the peaks of $r(t)$. Fig. 4d-4f present the CPU time, number of active constraints $c^*$ at the optimal solution $\theta^*$ and number of iterations performed $m^*$ of imuQP and the state-of-the-art active-set QP solvers qpOASES, quadprog and ODYS. Hildreth and mpcActiveSet-Solver are excluded as they are too slow. In Fig. 4e, the size of $\theta \in \mathbb{R}^{Nn_u}$, i.e. $Nn_u = 162$, is shown as a dashed line. It can be seen that $c^* \leq Nn_u = 162$ (or $n_\mathbb{A} \leq n = 162$ in [8, ch. 4.6.1, p. 62]), so that the system of linear equations (active constraints that hold with equality) is not overdetermined.

TABLE V
ACCURACY MEASURES: 1. *stationarity*; 2. *primal feasibility*; 3. *dual feasibility* AND *complementary slackness*; 4. *error w.r.t. quadprog*. THE ABBREVIATIONS ARE: ACC.: ACCORDING.

| | |
|---|---|
| 1 | $\sqrt{\sum_{i=1}^{Nn_u} (\mathcal{V}_i(\theta^*, \lambda^*))^2} = 0$ acc. to (20), $\mathcal{V}(\theta, \lambda) = E\theta + F + M^T\lambda$ |
| 2 | $\sqrt{\sum_{\mathcal{W}_i > 0} (\mathcal{W}_i(\theta^*))^2} = 0$ acc. to (15), $\mathcal{W}(\theta) = M\theta - \gamma$ |
| 3 | $\sqrt{\sum_{\lambda_i^* < 0} (\lambda_i^*)^2} = 0$ and $|\lambda^*|^T |M\theta^* - \gamma| = 0$ acc. to (15) |
| 4 | $\sqrt{\sum_{i=1}^{Nn_u} (\theta_i^* - \theta_{i,\text{quadprog}}^*)^2}$, the closer to 0, the closer $\theta^*$ is to $\theta_{\text{quadprog}}^*$ |

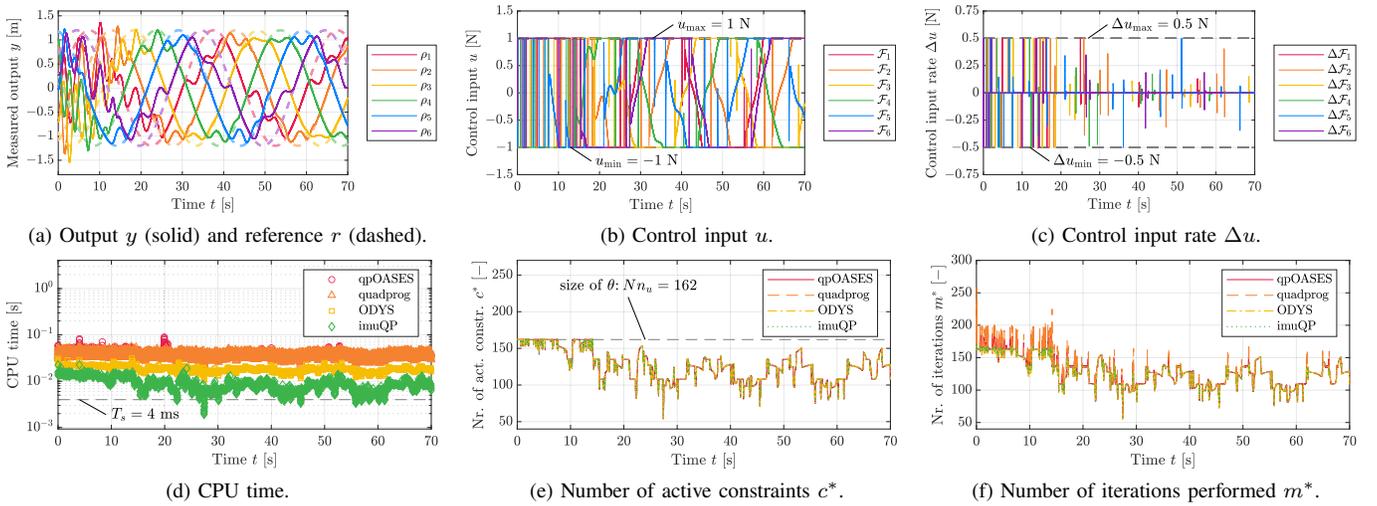

Fig. 4. Simulating a chain of six connected masses for $N = 27$ with single initial condition, simulated time is 70 s with 17500 samples. Simulated trajectories with imuQP as QP solver in (a) to (c), where $\rho_i$, $\mathcal{F}_i$ and $\Delta\mathcal{F}_i$ have the same color as mass $m_i$ for $i = 1,\ldots,6$ in Fig. 3. Comparing speeds between imuQP and state-of-the-art active-set QP solvers in (d) to (f). The abbreviations are: nr.: number; act.: active; constr.: constraints.

Speed and accuracy for each QP solver are presented in the top half of Table VII, where the lowest value in each column is highlighted in **bold**. The average, maximum and minimum CPU time (columns 3-5) are calculated with $k_a = 33$ and $k_b = 17500$ using (60), where the first few samples are skipped which usually show undesired spikes due to initialization, see Fig. 4d. The total CPU time and average values (columns 6-12) are calculated over all samples with $k_a = 1$ and $k_b = 17500$ using (60). The dual feasibility is 0 for all QP solvers and is therefore not shown. In Fig. 4d and Table VII, it can be seen that imuQP is the fastest of the four QP solvers and the only solver that gets below $T_s = 4$ ms — its average CPU time has improved by 51.37% compared with ODYS — while it performs around the same number of iterations and finds the same number of active constraints as the other QP solvers. Even though imuQP is not the most accurate solver, its accuracy is in the same order of magnitude as qpOASES'.

The 1-year free-trial version of the ODYS QP solver [1] has a maximum limit of 300 decision variables and a total of 1000 affine (inequality or equality) constraints. Therefore, we chose $N = 27$ with $p = 984$ from (13) and $\theta \in \mathbb{R}^{162}$.

### B. $N = 27$, many initial conditions

The same simulation settings from Section V-A apply, i.e. 70 s, 17500 samples, $T_s = 4$ ms, but now 35 simulations were run for different initial positions $\rho_i(0)$ using the following seeds: $\{[1,5], [7,22], 26, 28, 29, [31,35], 37, 38, [40,43]\}$. Some seeds are left out where qpOASES gives an error that the maximum number of working set calculations were performed.

The measurements in terms of speed and accuracy, averaged over the 35 simulations, are summarized in the bottom half of Table VII with the same $k_a$ and $k_b$ using (60) as in Section V-A. Again, the dual feasibility is 0 for all QP solvers and is therefore not shown. The lowest value in each column is again highlighted in **bold**, which almost has the same pattern as the single simulation in Section V-A and the same observations can be made regarding speed and accuracy. The average CPU time of imuQP, averaged over the 35 simulations, has improved by 51.93% compared with ODYS.

### C. Increasing N, same initial condition

For increasing $N$, 35 simulations are run for 15 s each — or 3750 samples with $T_s = 4$ ms — each with the same seed of 5 for the initial position $\rho_i(0)$ as in Section V-A. For Hildreth we set $\bar{m} = 1 \cdot 10^9$ and the tolerance $\delta$ is chosen differently for each $N$, for the same numerical precision reasons mentioned in Lemma 3.1, as follows: $\delta = 1 \cdot 10^{-29}$ for $N \in [1,5]$, $\delta = 1 \cdot 10^{-28}$ for $N = 6$, $\delta = 1 \cdot 10^{-27}$ for $N = 7$, $\delta = 1 \cdot 10^{-26}$ for $N \in [8,16]$ and $\delta = 1 \cdot 10^{-25}$ for $N \in [17,21]$.

Hildreth was simulated for $N \in [1,21]$, because it is too slow for larger $N$. For the same reason, mpcActiveSetSolver was run for $N \in [1, 25]$. ODYS was performed for $N \in [1, 27]$, because of the limitations in the 1-year free-trial version mentioned before. The QP solvers qpOASES, quadprog and imuQP were run for $N \in [1, 38]$, but with a few remarks. First, $N \in \{33, 35, 37\}$ were excluded for qpOASES, quadprog and imuQP, because quadprog's accuracy was very low. Second, $N \in \{11, [15, 20], 23, 24, 26, 34, 35\}$ were excluded for qpOASES, because it gave the same error that the maximum number of working set calculations were performed.

TABLE VI
QP SOLVER VERSIONS AND SETTINGS. DEFAULT USED IF UNSPECIFIED.

| QP solver | Version | Settings |
|---|---|---|
| HILDRETH | From [11], ch. 2.4.5, p. 67-68] | `lambda = zeros(height_M,1);`<br>`[theta,m,lambda] = hildreth(E,F,M,gamma,...`<br>`    lambda);` |
| MPCACTIVE-SETSOLVER | MATLAB R2021b | `iA0 = false(height_M,1);`<br>`Aeq = zeros(0,size_theta); beq = zeros(0,1);`<br>`options = mpcActiveSetOptions;`<br>`[theta,m,~,ineq] = mpcActiveSetSolver(E,F,...`<br>`    M,gamma,Aeq,beq,iA0,options);`<br>`lambda = ineq.ineqlin;` |
| QPOASES | 3.2.1 from [10] | `myOptions = qpOASES_options('fast');`<br>`myOptions.terminationTolerance = 1e-23;`<br>`[theta,~,~,m,ineq,~] = qpOASES(E,F,M,[],...`<br>`    [],[],gamma,myOptions);`<br>`lambda = -ineq(end - height_M + 1:end);` |
| QUADPROG | MATLAB R2021b | `theta_start = zeros(size_theta,1);`<br>`options = optimoptions('quadprog',...`<br>`    'Algorithm','active-set','Display','off');`<br>`[theta,~,~,output,ineq] = quadprog(E,F,M,...`<br>`    gamma,[],[],[],[],theta_start,options);`<br>`m = output.iterations; lambda = ineq.ineqlin;` |
| ODYS | 4.1.8 from [1] | `opt.returny = uint8(1);`<br>`[theta,~,~,ineq,info] = odysqpd(E,F,M,[],...`<br>`    gamma,[],[],[],[],opt);`<br>`m = info(1); lambda = ineq.ineq;` |
| IMUQP | 1 | `[theta,m,c,act,lambda_act] = imuQP(E_LU,F,...`<br>`    M,gamma,H); lambda = zeros(height_M,1);`<br>`lambda(act(1:c)) = lambda_act(1:c);` |

TABLE VII
SPEED AND ACCURACY OF IMUQP AND STATE-OF-THE-ART ACTIVE-SET QP SOLVERS FOR $N = 27$, SIMULATED TIME IS 70 S WITH 17500 SAMPLES: 1 SIMULATION (TOP) AND AVERAGED OVER 35 SIMULATIONS EACH WITH A DIFFERENT INITIAL CONDITION (BOTTOM). THE ABBREVIATIONS AND SYMBOLS ARE: AVG.: AVERAGE; MAX.: MAXIMUM; MIN.: MINIMUM; NR.: NUMBER; SIM.: SIMULATION; M: MINUTES; S: SECONDS.

|  | QP solver | Avg. CPU time [s] | Max. CPU time [s] | Min. CPU time [s] | Total CPU time | Avg. nr. of active constraints [−] | Avg. nr. of iterations [−] | Average stationarity [−] | Average primal feasibility [−] | Average complementary slackness [−] | Average error w.r.t. quadprog [N] |
|---|---|---|---|---|---|---|---|---|---|---|---|
| 1 sim. | QPOASES | 0.0404 | 0.0899 | 0.0195 | 11 M 47 S | 126.5562 | **127.3074** | $1.6866 \cdot 10^{-13}$ | $6.4546 \cdot 10^{-12}$ | $1.1877 \cdot 10^{-10}$ | **$7.4312 \cdot 10^{-11}$** |
|  | QUADPROG | 0.0340 | 0.0583 | 0.0194 | 9 M 55 S | 126.5562 | 129.5932 | **$3.6149 \cdot 10^{-14}$** | **$3.1327 \cdot 10^{-16}$** | **$3.5685 \cdot 10^{-16}$** | — |
|  | ODYS | 0.0183 | 0.0323 | 0.0099 | 5 M 20 S | 126.5562 | 127.6078 | $1.0248 \cdot 10^{-12}$ | $6.6111 \cdot 10^{-13}$ | $2.4703 \cdot 10^{-11}$ | $1.3391 \cdot 10^{-10}$ |
|  | IMUQP | **0.0089** | **0.0225** | **0.0019** | **2 M 35 S** | 126.5562 | 127.6571 | $8.7136 \cdot 10^{-14}$ | $8.3868 \cdot 10^{-12}$ | $5.8611 \cdot 10^{-10}$ | $4.9164 \cdot 10^{-10}$ |
| 35 sims. | QPOASES | 0.0398 | 0.0765 | 0.0233 | 11 M 37 S | 125.7436 | **126.4285** | $1.6041 \cdot 10^{-13}$ | $6.4448 \cdot 10^{-12}$ | $1.1371 \cdot 10^{-10}$ | $5.9147 \cdot 10^{-10}$ |
|  | QUADPROG | 0.0340 | 0.0620 | 0.0224 | 9 M 56 S | 125.7436 | 128.7010 | **$3.4997 \cdot 10^{-14}$** | **$3.1854 \cdot 10^{-16}$** | **$3.2910 \cdot 10^{-16}$** | — |
|  | ODYS | 0.0181 | 0.0302 | 0.0110 | 5 M 18 S | 125.7436 | 126.7883 | $1.0124 \cdot 10^{-12}$ | $5.9849 \cdot 10^{-13}$ | $2.1309 \cdot 10^{-11}$ | **$3.0340 \cdot 10^{-10}$** |
|  | IMUQP | **0.0087** | **0.0228** | **0.0028** | **2 M 33 S** | 125.7436 | 126.8309 | $8.4655 \cdot 10^{-14}$ | $7.8400 \cdot 10^{-12}$ | $5.2177 \cdot 10^{-10}$ | $1.3036 \cdot 10^{-8}$ |

The speed and accuracy are plotted as functions of $N$ in Fig. 5, where each data point corresponds to a single simulation. The idea of plotting primal feasibility as a function of $N$ was taken from [15]. In Fig. 5a, the average, maximum and minimum CPU time are calculated with $k_a = 33$ and $k_b = 3750$ using (60) for the same reasons mentioned already. The total CPU time and other average values − in Fig. 5b-5h − are calculated over all samples with $k_a = 1$ and $k_b = 3750$ using (60). The dual feasibility is 0 for each $N$ and all QP solvers and is therefore not shown. It can be seen that, as $N$ increases, the CPU time of imuQP increases with a slower rate than the other QP solvers. The average stationarity of imuQP lies between that of qpOASES and quadprog for increasing $N$. However, the average primal feasibility and average complementary slackness get worse faster than the other QP solvers, but still remain within the same order of magnitude. Again, imuQP finds the same number of active constraints as the other QP solvers. In Fig. 5a, it can be seen that the maximum $N$ before the maximum CPU time exceeds $T_s = 4$ ms is $N = 2$ for Hildreth, $N = 4$ for mpcActiveSetSolver, $N = 6$ for quadprog, $N = 8$ for qpOASES, $N = 13$ for ODYS and $N = 14$ for imuQP. In Fig. 5c, the theoretical maximum number of active constraints $w = \frac{p}{2}$ in (13) and the size of $\theta \in \mathbb{R}^{Nn_u}$, i.e. $Nn_u$, are shown as dashed lines as a function of $N$. Except for Hildreth due to low accuracy in $\lambda^*$, for the other QP solvers it can be seen that $c^* \leq Nn_u$ as discussed in Section V-A.

*Remark 8.* For a real-time implementation and experimental validation of an earlier version of imuQP − which excludes infeasibility detection and linear dependence check of constraints − on an Arduino Due for integral-action MPC on an educational magnetic ball levitation setup with sampling period $T_s = 5$ ms, we refer to [48] due to space limitations.

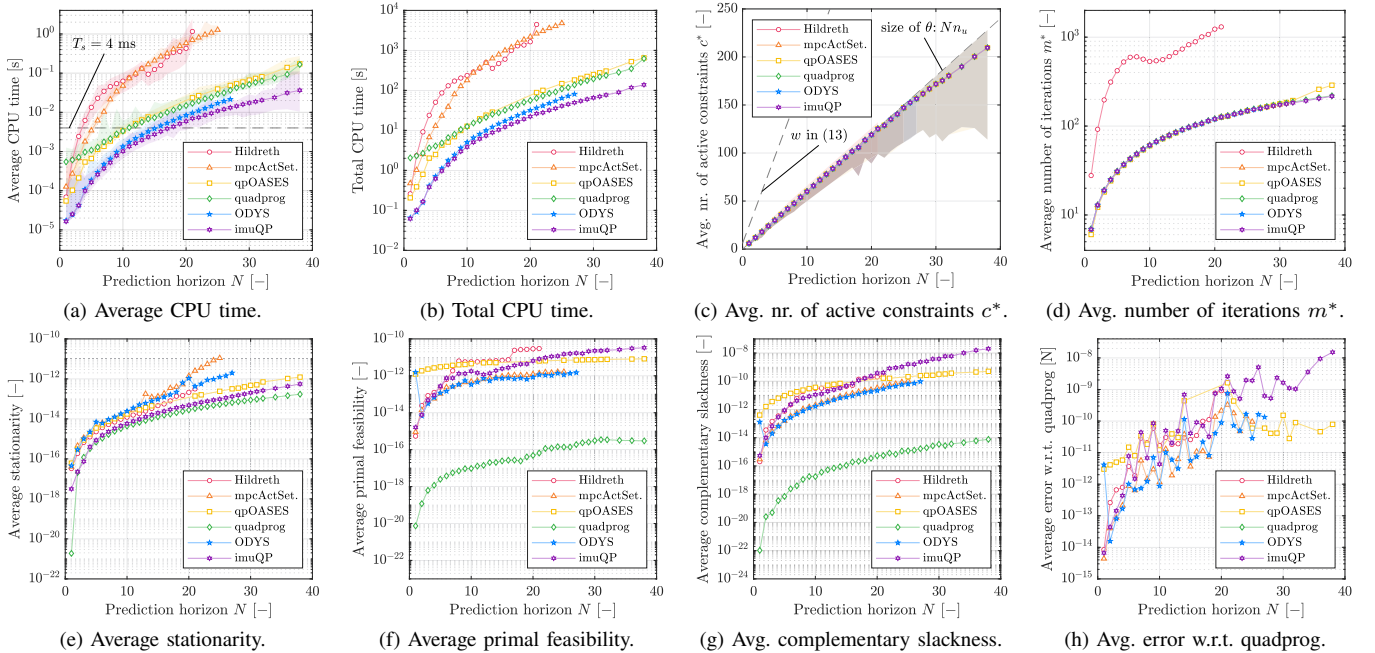

Fig. 5. Speed and accuracy of imuQP and state-of-the-art active-set QP solvers in simulation for increasing $N$ with the same initial conditions, simulated time is 15 s with 3750 samples. In (a) and (c), the min. and max. values are shown as transparent colors. Note that in (c) and (d) some data points are overlapping. The abbreviations are: mpcActSet.: mpcActiveSetSolver; avg.: average; nr.: number; min.: minimum; max.: maximum.

(a) Average CPU time.  (b) Total CPU time.  (c) Avg. nr. of active constraints $c^*$.  (d) Avg. number of iterations $m^*$.
(e) Average stationarity.  (f) Average primal feasibility.  (g) Avg. complementary slackness.  (h) Avg. error w.r.t. quadprog.

## VI. CONCLUSIONS

THIS paper presented a new fast active-set QP solver, called imuQP. It is based on KKT conditions and inverse matrix updates, which compute the Lagrange multiplier $\lambda_{\mathbb{A}}$ cheaply and very fast for each update of the active set. The algorithm and convergence analysis of imuQP including infeasibility detection were presented and its computational complexity and memory footprint were compared with qpOASES, which is a well-known active-set QP solver.

Speed and accuracy of imuQP were compared with state-of-the-art active-set QP solvers through MATLAB simulations, where integral-action MPC was performed on a chain of six spring-connected masses with $T_s = 4$ ms. For $N = 27$, It is shown that the average CPU time of imuQP has improved by at least 51% compared with ODYS. The accuracy of imuQP is of the same order of magnitude as qpOASES'. As $N$ increases, the CPU time of imuQP increases slower while its accuracy gets worse faster but still remains within the same order of magnitude, compared with other QP solvers. Therefore, imuQP is competitive with state-of-the-art active-set QP solvers and is suitable for real-time MPC for fast systems with small $T_s$.

## APPENDIX A: COMPACT MPC CONSTRAINT FORMULATION

We can write the constraints on $\Delta u_{i|k}$ and $y_{i|k}$ in (8) into the general form for $i = 0, \ldots, N-1$

$$Z_i \Delta u_{i|k} + T_i y_{i|k} \le b_i, \quad T_N y_{N|k} \le b_N, \quad (A.1)$$

where

$$Z_i = \begin{bmatrix} -I_{n_u} \\ I_{n_u} \\ 0 \\ 0 \end{bmatrix}, \quad T_i = \begin{bmatrix} 0 \\ 0 \\ -I_{n_y} \\ I_{n_y} \end{bmatrix}, \quad b_i = \begin{bmatrix} -\Delta u_{\min} \\ \Delta u_{\max} \\ -y_{\min} \\ y_{\max} \end{bmatrix},$$

$$T_N = \begin{bmatrix} -I_{n_y} \\ I_{n_y} \end{bmatrix}, \quad b_N = \begin{bmatrix} -y_{\min} \\ y_{\max} \end{bmatrix},$$

and the zero matrix 0 is of appropriate dimensions. Using $y(k) = Cx(k)$, (A.1) can be rewritten in compact form as

$$\mathcal{Z} \Delta U_k + \mathcal{D} x(k) + \mathcal{T} Y_k \le \mathcal{C}, \quad (A.2)$$

where

$$\mathcal{Z} = \begin{bmatrix} Z_0 & \cdots & 0 \\ \vdots & \ddots & \vdots \\ 0 & \cdots & Z_{N-1} \\ 0 & \cdots & 0 \end{bmatrix}, \quad \mathcal{D} = \begin{bmatrix} T_0 C \\ 0 \\ \vdots \\ 0 \end{bmatrix},$$

$$\mathcal{T} = \begin{bmatrix} 0 & \cdots & 0 \\ T_1 & \cdots & 0 \\ \vdots & \ddots & \vdots \\ 0 & \cdots & T_N \end{bmatrix}, \quad \mathcal{C} = \begin{bmatrix} b_0 \\ b_1 \\ \vdots \\ b_N \end{bmatrix}.$$

The constraints on $u_{i|k}$ in (8) can be compactly written as

$$\underbrace{\begin{bmatrix} u_{\min} \\ \vdots \\ u_{\min} \end{bmatrix}}_{U_{\min}} \le \underbrace{\begin{bmatrix} u_{0|k} \\ \vdots \\ u_{N-1|k} \end{bmatrix}}_{U_k} \le \underbrace{\begin{bmatrix} u_{\max} \\ \vdots \\ u_{\max} \end{bmatrix}}_{U_{\max}}, \quad (A.3)$$

where $U_k \in \mathbb{R}^{Nn_u}, U_{\min} \in \mathbb{R}^{Nn_u}, U_{\max} \in \mathbb{R}^{Nn_u}$ and the relation between $U_k$ and $\Delta U_k$ is

$$U_k = \underbrace{\begin{bmatrix} I_{n_u} \\ \vdots \\ I_{n_u} \end{bmatrix}}_{\mathcal{S}} u(k-1) + \underbrace{\begin{bmatrix} I_{n_u} & \cdots & 0 \\ \vdots & \ddots & \vdots \\ I_{n_u} & \cdots & I_{n_u} \end{bmatrix}}_{\mathcal{K}} \Delta U_k. \quad (A.4)$$

Substituting $Y_k$ in (4) into (A.2) and substituting $U_k$ in (A.4) into (A.3), we obtain the affine inequality constraints

$$\mathcal{L} \Delta U_k \le d + W x(k) + V u(k-1), \quad (A.5)$$

where

$$\mathcal{L} = \begin{bmatrix} \mathcal{T} \Gamma + \mathcal{Z} \\ -\mathcal{K} \\ \mathcal{K} \end{bmatrix}, \quad d = \begin{bmatrix} \mathcal{C} \\ -U_{\min} \\ U_{\max} \end{bmatrix},$$

$$W = \begin{bmatrix} -\mathcal{D} - \mathcal{T} \Phi \\ 0 \\ 0 \end{bmatrix}, \quad V = \begin{bmatrix} 0 \\ \mathcal{S} \\ -\mathcal{S} \end{bmatrix},$$

and the zero matrix 0 is of appropriate dimensions, $\mathcal{C}$ is defined in (A.2) and $\mathcal{S}, \mathcal{K}$ are defined in (A.4).

## ACKNOWLEDGEMENT

The authors would like to thank dr. Thành Mừng Lâm for suggesting the sinusoidal reference signal $r_i(t)$ in (59) and Fig. 4a and MSc Yuk Hang Yuen for his contribution in improving the Hildreth algorithm and suggesting the use of MEX files. The YouTube lectures by dr. Richard G. Brown from the Massey University of New Zealand inspired the design of imuQP and the proof of Lemma 2.3.(i) together with Fig. 1.

## REFERENCES

[1] G. Cimini, A. Bemporad, and D. Bernardini, "ODYS QP Solver," *ODYS S.r.l.*, 2017. [Online]. Available: https://www.odys.it/qp-solver-for-embedded-optimization
[2] B. Stellato, G. Banjac, P. Goulart, A. Bemporad, and S. P. Boyd, "OSQP: an operator splitting solver for quadratic programs," *Mathematical Programming Computation*, vol. 12, no. 4, pp. 637–672, 2020.
[3] J. V. Frasch, M. Vukov, H. J. Ferreau, and M. Diehl, "A dual Newton strategy for the efficient solution of sparse quadratic programs arising in SQP-based nonlinear MPC," *Optimization Online*, 2013.
[4] ——, "A new quadratic programming strategy for efficient sparsity exploitation in SQP-based nonlinear MPC and MHE," *IFAC Proceedings Volumes*, vol. 47, no. 3, pp. 2945–2950, 2014.
[5] J. V. Frasch, S. Sager, and M. Diehl, "A parallel quadratic programming method for dynamic optimization problems," *Mathematical Programming Computation*, vol. 7, no. 3, pp. 289–329, 2015.
[6] J. V. Frasch, H. J. Ferreau, M. Diehl, and M. Vukov, "qpDUNES – a DUal NEwton Strategy for convex quadratic programming," *GitHub*, 2016. [Online]. Available: https://github.com/jfrasch/qpDUNES
[7] H. J. Ferreau, "Model Predictive Control Algorithms for Applications with Millisecond Timescales," *PhD thesis, Department of Electrical Engineering, Katholieke Universiteit Leuven, Leuven, Belgium*, 2011.
[8] ——, "An Online Active Set Strategy for Fast Solution of Parametric Quadratic Programs with Applications to Predictive Engine Control," *MSc thesis, Faculty of Mathematics and Computer Science, Heidelberg University, Heidelberg, Germany*, 2006.
[9] H. J. Ferreau, C. Kirches, A. Potschka, H. G. Bock, and M. Diehl, "qpOASES: a parametric active-set algorithm for quadratic programming," *Math. Program. Computation*, vol. 6, no. 4, pp. 327–363, 2014.
[10] H. J. Ferreau *et al.*, "qpOASES User's Manual, Version 3.2 (April 2017)," *ABB Corporate Research, Switzerland*, 2017. [Online]. Available: https://www.coin-or.org/qpOASES/doc/3.2/manual.pdf
[11] L. Wang, "Model Predictive Control System Design and Implementation Using MATLAB®," *In M. J. Grimble and M. A. Johnson (Eds.), Advances in Industrial Control, Springer*, 2009.

[12] C. Hildreth, "A quadratic programming procedure," *Naval Research Logistics Quarterly*, vol. 4, no. 1, pp. 79–85, 1957.
[13] A. N. Iusem and A. R. De Pierro, "On the convergence properties of Hildreth's quadratic programming algorithm," *Mathematical Programming*, vol. 47, no. 1, pp. 37–51, 1990.
[14] A. Lent and Y. Censor, "Extensions of Hildreth's Row-action Method for Quadratic Programming," *SIAM Journal on Control and Optimization*, vol. 18, no. 4, pp. 444–454, 1980.
[15] Y. J. J. Heuts, G. P. Padilla, and M. C. F. Donkers, "An Adaptive Restart Heavy-Ball Projected Primal-Dual Method for Solving Constrained Linear Quadratic Optimal Control Problems," *2021 60$^{th}$ IEEE Conference on Decision and Control (CDC)*, pp. 6722–6727, 2021.
[16] C. Schmid and L. T. Biegler, "Quadratic Programming Methods for Reduced Hessian SQP," *Computers & Chemical Engineering*, vol. 18, no. 9, pp. 817–832, 1994.
[17] D. Piga, S. Formentin, and A. Bemporad, "Direct Data-Driven Control of Constrained Systems," *IEEE Transactions on Control Systems Technology*, vol. 26, no. 4, pp. 1422–1429, 2018.
[18] P. Chalupa, J. Novák, and M. Malý, "Modelling and model predictive control of magnetic levitation laboratory plant," *Proc. 31$^{st}$ European Conf. on Modelling and Simulation (ECMS 2017)*, pp. 367–373, 2017.
[19] P. Zometa, M. Kögel, T. Faulwasser, and R. Findeisen, "Implementation Aspects of Model Predictive Control for Embedded Systems," *2012 American Control Conference (ACC)*, pp. 1205–1210, 2012.
[20] R. Kouki, H. Salhi, and F. Bouani, "Application of Model Predictive Control for a thermal process using STM32 Microcontroller," *2017 International Conference on Control, Automation and Diagnosis (ICCAD'17)*, pp. 146–151, 2017.
[21] pronenewbits, "Arduino Unconstrained MPC Library," *GitHub*, 2020. [Online]. Available: https://github.com/pronenewbits/Arduino_Unconstrained_MPC_Library
[22] M. Gulan, G. Takács, N. A. Nguyen, S. Olaru, P. Rodríguez-Ayerbe, and B. Rohal'-Ilkiv, "Efficient Embedded Model Predictive Vibration Control via Convex Lifting," *IEEE Transactions on Control Systems Technology*, vol. 27, no. 1, pp. 48–62, 2019.
[23] APMonitor, "Process Control Temperature Lab," *GitHub*, 2021. [Online]. Available: https://github.com/APMonitor/arduino
[24] S. Zavitsanou, A. Chakrabarty, E. Dassau, and F. J. Doyle III, "Embedded Control in Wearable Medical Devices: Application to the Artificial Pancreas," *Processes*, vol. 4, no. 4, pp. 1–29, 2016.
[25] B. J. T. Binder, D. K. M. Kufoalor, and T. A. Johansen, "Scalability of QP solvers for Embedded Model Predictive Control Applied to a Subsea Petroleum Production System," *2015 IEEE Conference on Control Applications (CCA)*, pp. 1173–1178, 2015.
[26] D. K. M. Kufoalor, B. J. T. Binder, H. J. Ferreau, L. Imsland, T. A. Johansen, and M. Diehl, "Automatic Deployment of Industrial Embedded Model Predictive Control using qpOASES," *2015 European Control Conference (ECC)*, pp. 2601–2608, 2015.
[27] B. Huyck, L. Callebaut, F. Logist, H. J. Ferreau, M. Diehl, J. De Brabanter, J. Van Impe, and B. De Moor, "Implementation and Experimental Validation of Classic MPC on Programmable Logic Controllers," *2012 20$^{th}$ Mediterr. Conf. on Control & Autom. (MED)*, pp. 679–684, 2012.
[28] B. Huyck, H. J. Ferreau, M. Diehl, J. De Brabanter, J. F. M. Van Impe, B. De Moor, and F. Logist, "Towards Online Model Predictive Control on a Programmable Logic Controller: Practical Considerations," *Mathematical Problems in Engineering*, vol. 2012, pp. 1–20, 2012.
[29] B. Huyck, J. De Brabanter, B. De Moor, J. F. Van Impe, and F. Logist, "Online model predictive control of industrial processes using low level control hardware: A pilot-scale distillation column case study," *Control Engineering Practice*, vol. 28, pp. 34–48, 2014.
[30] C. Ibañez, C. Ocampo-Martinez, and B. Gonzalez, "Embedded optimization-based controllers for industrial processes," *2017 IEEE 3$^{rd}$ Colombian Conference on Automatic Control (CCAC)*, 2017.
[31] C. Ibáñez López, "Implementation of optimization-based controllers for industrial processes," *MSc thesis, Barcelona School of Industrial Engineering, Polytechnic University of Catalonia, Barcelona, Spain*, 2017.
[32] S. Adhau, S. Patil, D. Ingole, and D. Sonawane, "Implementation and Analysis of Nonlinear Model Predictive Controller on Embedded Systems for Real-Time Applications," *2019 18$^{th}$ European Control Conference (ECC)*, pp. 3359–3364, 2019.
[33] V. T. T. Lam, A. Sattar, L. Wang, and M. Lazar, "Fast Hildreth-based Model Predictive Control of Roll Angle for a Fixed-Wing UAV," *IFAC-PapersOnLine*, vol. 53, no. 2, pp. 5757–5763, 2020.
[34] Arduino, "Arduino Due," *Arduino Official Store*, 2024. [Online]. Available: https://store.arduino.cc/arduino-due
[35] D. G. Luenberger and Y. Ye, "Linear and Nonlinear Programming, Fourth Edition," *In C. C. Price (Ed.), International Series in Operations Research & Management Science, Springer*, vol. 228, 2016.
[36] J. Nocedal and S. J. Wright, "Numerical Optimization, Second Edition," *In T. V. Mikosch, S. M. Robinson and S. I. Resnick (Eds.), Springer Series in Operations Research, Springer*, 2006.
[37] pronenewbits, "Arduino Constrained MPC Library," *GitHub*, 2020. [Online]. Available: https://github.com/pronenewbits/Arduino_Constrained_MPC_Library
[38] Y. Cao, "Update Inverse Matrix: Version 1.0.0.0," *MATLAB Central File Exchange*, 2008. [Online]. Available: https://mathworks.com/matlabcentral/fileexchange/18063-update-inverse-matrix
[39] W. W. Hager, "Updating the Inverse of a Matrix," *SIAM Review*, vol. 31, no. 2, pp. 221–239, 1989.
[40] M. E. Khan, "Updating Inverse of a Matrix When a Column is Added/Removed," *Technical report, Department of Computer Science, University of British Columbia, Vancouver, Canada*, 2008. [Online]. Available: https://emtiyaz.github.io/Writings/OneColInv.pdf
[41] J. Gallier, "The Schur Complement and Symmetric Positive Semidefinite (and Definite) Matrices," *Technical report, Computer and Information Science Dept., University of Pennsylvania, Philadelphia, USA*, 2019. [Online]. Available: http://www.cis.upenn.edu/~jean/schur-comp.pdf
[42] R. A. Horn and C. R. Johnson, "Matrix Analysis, Second Edition," *Cambridge University Press*, 2013.
[43] S. Boyd and L. Vandenberghe, "Convex Optimization," *Cambridge University Press*, 2004.
[44] J. Matoušek and B. Gärtner, "Understanding and Using Linear Programming," *In Universitext, Springer*, 2007.
[45] G. H. Golub and C. F. Van Loan, "Matrix Computations, Third Edition," *Johns Hopkins University Press*, 1996.
[46] M. Kögel, P. Zometa, and R. Findeisen, "On Tailored Model Predictive Control for Low Cost Embedded Systems with Memory and Computational Power Constraints," *Technical report, Institute for Automation Engineering, Otto von Guericke University Magdeburg, Magdeburg, Germany*, 2012.
[47] D. Arnström, A. Bemporad, and D. Axehill, "A Dual Active-Set Solver for Embedded Quadratic Programming Using Recursive $LDL^T$ Updates," *IEEE Transactions on Automatic Control*, vol. 67, no. 8, pp. 4362–4369, 2022.
[48] V. T. T. Lâm, "QP-MPC solvers for real-time control of fast nonlinear systems using inexpensive microcontrollers," *MSc thesis, Department of Electrical Engineering, Eindhoven University of Technology, Eindhoven, the Netherlands*, 2020. [Online]. Available: https://pure.tue.nl/ws/portalfiles/portal/261832470/0857216_Lam_VTT_Thesis_paperFinalV6.pdf

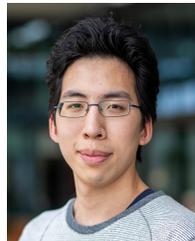

**Victor Trường Thịnh Lâm** was born on March 14$^{th}$, 1995 in Dordrecht, the Netherlands and obtained his BSc degree in January 2017 and his MSc degree in the Control Systems group in July 2020, both at the Department of Electrical Engineering, Eindhoven University of Technology, Eindhoven, the Netherlands. Since February 2022, he is a PhD student in the Electromechanics and Power Electronics group at the same department. His research interest is in model predictive control for mechatronic systems.

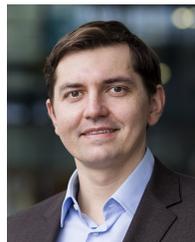

**Mircea Lazar** was born on March 4$^{th}$, 1978 in Iași, Romania and is an associate professor in Constrained Control of Complex Systems at the Department of Electrical Engineering, Eindhoven University of Technology, Eindhoven, the Netherlands. His research interests cover stability theory, Lyapunov functions, distributed control, and constrained control of non-linear and hybrid systems, including model predictive control. His research is driven by control problems in power systems, power electronics, high-precision mechatronics, automotive and biological systems. Dr. Lazar received the European Embedded Control Institute PhD Award in 2007 for his PhD dissertation. He received a VENI personal grant from the Dutch Research Council (NWO) in 2008 and he supervised 10 PhD researchers (2 received the cum laude distinction). Dr. Lazar chaired the 4$^{th}$ IFAC Conference on Non-Linear Model Predictive Control in 2012. He is an active member of the IFAC Technical Committees 1.3 Discrete Event and Hybrid Systems, 2.3 Non-Linear Control Systems and an associate editor of IEEE Transactions on Automatic Control.